\documentclass{article}

\begin{document}

\title{Atiyah's work on holomorphic vector bundles and gauge theories}
\author{Simon Donaldson}
\date{\today}
\maketitle


\tableofcontents
\newcommand{\inter}{{\rm int}}
\newcommand{\bC}{{\bf C}}
\newcommand{\bP}{{\bf P}}
\newcommand{\bR}{{\bf R}}
\newcommand{\bQ}{{\bf Q}}
\newcommand{\bH}{{\bf H}}
\newcommand{\cO}{{\cal O}}
\newcommand{\rHom}{{\rm Hom}}
\newcommand{\db}{\overline{\partial}}
\newcommand{\Vect}{{\rm Vect}}
\newcommand{\cL}{{\cal L}}
\newcommand{\bZ}{{\bf Z}}
\newcommand{\Riem}{{\rm Riem}}
\newcommand{\ric}{{\rm Ricci}}
\newcommand{\uE}{\underline{E}}
\newcommand{\uU}{\underline{U}}
\newcommand{\uV}{\underline{V}}
\newcommand{\uW}{\underline{W}}
\newcommand{\uH}{\underline{H}}
\newcommand{\uM}{\underline{M}}
\newcommand{\oM}{\overline{M}}
\newcommand{\cV}{{\cal V}}
\newcommand{\cM}{{\cal M}}
\newcommand{\cA}{{\cal A}}
\newcommand{\cG}{{\cal G}}
\newcommand{\cB}{{\cal B}}
\newcommand{\rad}{{\rm ad}}
\newcommand{\tM}{\tilde{M}}
\newcommand{\tS}{\tilde{S}}
\newcommand{\tX}{\tilde{X}}
\newcommand{\ubC}{\underline{{\bf C}}}
\newcommand{\bT}{{\bf T}}
\newcommand{\cE}{{\cal E}}
\newcommand{\Jac}{{\rm Jac}}
\newcommand{\unabla}{\underline{\nabla}}
\newcommand{\bE}{{\bf E}}
\newcommand{\uLambda}{\underline{\Lambda}}
\newcommand{\frakg}{{\rm g}}
\newcommand{\bF}{{\bf F}}
\newtheorem{defn}{Definition}
\newtheorem{prop}{Proposition}
\newtheorem{thm}{Theorem}
\section{Early work}

Atiyah began his research career as an algebraic geometer, writing a series of papers in the period 1952-1958, mostly about bundles over algebraic varieties.  Many of the themes in these papers reappear 20 years later when Atiyah's interests turned to  Yang-Mills theories. In the commentary on Volume 1 of his collected works Atiyah wrote: {\it I was fascinated by classical projective geometry\dots at the same time great things were happening in France and I was an avid  reader of the {\em Comptes Rendus}, following the developments in sheaf theory}. The language and style of these early papers of Atiyah are an interesting
mixture of classical and modern.

\subsection{Extensions and the Atiyah class}

{\it Extensions} of vector bundles form a theme running through much of Atiyah's work. Let $X$ be a complex manifold. A holomorphic vector bundle $\pi:E\rightarrow X$ is a complex vector bundle in the ordinary sense such that the total space is a complex manifold, and the projection and structure maps are all holomorphic. Alternatively the bundle is defined by a system of holomorphic transition functions 
$$  g_{\alpha\beta}: U_{\alpha}\cap U_{\beta}\rightarrow GL(r,\bC), $$ with respect to an open cover $X= \bigcup U_{\alpha}$. The definition of a subbundle $F\subset E$ is the obvious one and we get a holomorphic quotient $Q=E/F$ and an exact sequence:
\begin{equation}   0\rightarrow F\stackrel{i}{\rightarrow} E\stackrel{p}{\rightarrow} Q\rightarrow 0. \end{equation}
The new feature, compared with the theory of  topological or $C^{\infty}$ bundles, is that this exact sequence need not split: it is not usually true that $E$ is isomorphic to a direct sum $F\oplus Q$. We define an equivalence relation on short exact sequences like (1) by saying that two are equivalent if there is a commutative diagram:

\

\

\setlength{\unitlength}{.2pt}

\begin{picture}(1000,320)(-300,0)
\put(100,30){$0$}
\put(300,30){$F$}
\put(500,30){$E$}
\put(700,30){$Q$}
\put(900, 30){$0$}
\put(100,270){$0$}
\put(300,270){$F$}
\put(500,270){$E'$}
\put(700,270){$Q$}
\put(900,270){$0$}
\put(140,50){\vector(1,0){160}}
\put(340,50){\vector(1,0){160}}
\put(540,50){\vector(1,0){160}}
\put(740,50){\vector(1,0){160}}
\put(140,290){\vector(1,0){160}}
\put(340,290){\vector(1,0){160}}
\put(540,290){\vector(1,0){160}}
\put(740,290){\vector(1,0){160}}
\put(300,80){\line(0,1){160}}
\put(315,80){\line(0,1){160}}
\put(520,240){\vector(0,-1){160}}
\put(700,80){\line(0,1){160}}
\put(715,80){\line(0,1){160}}
\end{picture}

Then the basic fact is:

\begin{prop}
Given holomorphic bundles $F,Q$ there is a one-to-one correspondence between equivalence classes of extensions (1) and the sheaf cohomology group $H^{1}({\rm Hom}(Q,F))$ such that the trivial extension defined by $F\oplus Q$ corresponds to $0\in H^{1}({\rm
Hom}(Q,F))$.
\end{prop}

  In \cite{kn:MFA2} Atiyah writes that this result \lq\lq follows from the general theory of fibre bundles'' and refers to unpublished notes of Grothendieck. In the notation used above we do not distinguish between a holomorphic vector bundle and its sheaf of local holomorphic sections. There is a one-to-one correspondence between vector bundles and sheaves of locally free modules over the structure sheaf $\cO_{X}$ of  holomorphic functions on $X$ and a version of Proposition 1 holds for locally free sheaves of modules in general. To define the extension class from the axioms of sheaf cohomology we apply ${\rm Hom}(Q,\ )$ to the sequence to get 
$$   0\rightarrow {\rm Hom}(Q,F)\rightarrow {\rm Hom}(Q,E) \rightarrow {\rm Hom}(Q,Q)\rightarrow 0$$
and this has a long exact cohomology sequence with coboundary map
$$  \partial: H^{0}(\rHom(Q,Q))\rightarrow H^{1}(\rHom(Q,F)). $$
The extension class is $\partial 1_{Q}$ where $1_{Q}$ is the identity on each fibre of $Q$. 

We can see this more explicitly using either \v{C}ech or Dolbeault representations of the cohomology.  For the first we use the fact that the extension can be split locally, so there is an open cover $\{U_{\alpha}\}$ over which we have holomorphic splitting maps $$j_{\alpha}: Q\vert_{U_{\alpha}}\rightarrow E\vert_{U_{\alpha}}. $$ 
On an intersection $U_{\alpha\beta}=U_{\alpha}\cap U_{\beta}$ we have two splitting maps and 
$$ j_{\alpha}-j_{\beta}= i\circ T_{\alpha \beta}, $$
where $T_{\alpha \beta}: Q\vert_{U_{\alpha \beta}}\rightarrow F\vert_{U_{\alpha\beta}}$.
Thus $\{T_{\alpha\beta}\}$ is a  cochain in the \v{C}ech complex associated to the sheaf $\rHom(Q,F)$ and one finds that it is a cocycle and that its cohomology class is independent of choices. 

The Dolbeault description of the extension class is important in complex differential geometry. To set this up, we recall that a holomorphic vector bundle can be viewed as a $C^{\infty}$ bundle equipped with a $\db$-operator
$$  \db_{E}:\Omega^{0}(E)\rightarrow \Omega^{0,1}(E), $$
satisfying the Leibnitz rule $\db_{E}(fs)= (\db f) s + f \db_{E}s$. Here we are writing $\Omega^{0,1}(E)$ for the smooth $(0,1)$-forms with values in $E$. Using the Leibnitz rule one gets an extension to operators 
$$   \db_{E}:\Omega^{p,q}(E)\rightarrow \Omega^{p,q+1}(E)$$
with $\db_{E}^{2}=0$. Fixing $p$, the cohomology of the resulting complex gives the sheaf cohomology $H^{*}(E\otimes \Lambda^{p}T^{*}X)$.

Returning to the extension class, this time we use the fact that a sequence of $C^{\infty}$ bundles can be split.  So we can choose a global $C^{\infty}$ (not necessarily holomorphic) splitting map $j:Q\rightarrow E$ and a corresponding projection $\varpi:E\rightarrow F$. Then the composite $\varpi\circ \db_{E}\circ j$ is a map from $C^{\infty}$ sections of $Q$ to $(0,1)$-forms with values in  $F$. One finds that this commutes with multiplication by smooth functions so is given by a tensor $B\in \Omega^{0,1}(\rHom(Q,F))$ which gives a Dolbeault representative for the extension class. Said in another way, if we identify $E$ with $F\oplus Q$ as a $C^{\infty}$ bundle then $\db_{E}$ has a \lq\lq matrix'' representation:

\begin{equation} \db_{E}=\left(\begin{array}{cc} \db_{F}& B\\0&\db_{Q}\end{array}\right).\end{equation}

With this background in place we begin our review of some of Atiyah's applications of the theory. In \cite{kn:MFA2} he introduced  what is now called the {\it Atiyah class} of a bundle $E\rightarrow X$. This can be described in two ways. For the first, for a point $x\in X$ let $J(E)_{x}$ be the set of $1$-jets of sections of $E$ at $x$. By definition, such a $1$-jet is an equivalence class of local holomorphic sections where $s$ is equivalent to $s'$ if $s(x)=s'(x)$ and the derivative of $s-s'$ at $x$ vanishes. (Recall that the derivative of a section of a vector bundle is not intrinsically defined in general but it is defined at a zero of the section.) Then we have an exact sequence of vector spaces
$$   0\rightarrow T^{*}X_{x}\otimes E_{x}\rightarrow J(E)_{x}\rightarrow E_{x}\rightarrow 0.  $$
Letting $x$ vary we get an exact sequence of bundles over $X$
\begin{equation}  0\rightarrow T^{*}X\otimes E \rightarrow J(E)\rightarrow E\rightarrow
0,  \end{equation}
and the Atiyah class $\alpha(E)$ is the extension class in 
$$ H^{1}({\rm Hom}(E, E\otimes T^{*}X)= H^{1}({\rm End E}\otimes T^{*} X). $$

A more general  approach  uses principal bundles, and this was the point of view Atiyah took in \cite{kn:MFA2}. Let $G$ be a complex Lie group  and $\pi: P\rightarrow X$ a principal $G$-bundle. Then there is an exact sequence of vector bundles over $P$:
$$  0\rightarrow \underline{{\rm Lie}(G)}\rightarrow TP\rightarrow \pi^{*}TX\rightarrow 0$$
Here $\underline{{\rm Lie}(G)}$ is the trivial bundle with fibre the Lie algebra ${\rm Lie}(G)$ which can be identified with the \lq\lq vertical' tangent vectors along the fibres of $\pi$ using the group action. All three bundles in this sequence are equivariant with respect to the action of $G$ and this implies that the sequence is pulled back from an exact sequence of vector bundles over $X$
\begin{equation} 0\rightarrow {\rm ad}P\rightarrow \bT\rightarrow TX\rightarrow 0, \end{equation}
where ${\rm ad}P$ is the bundle over $X$ associated to the adjoint representation of $G$ on its Lie algebra. Then there is an extension class in $H^{1}(T^{*}X\otimes {\rm ad}P)$. 
In the case when $G=GL(r,\bC)$, so $P$ is the frame bundle of a vector bundle $E$, we have ${\rm ad}P= {\rm End} E$ and we get the same class as before. 
 The different approaches correspond to two ways of viewing {\it connections}. From the principal bundle point of view a connection is given by a splitting of the exact sequence (4), that is to say a $G$-invariant \lq\lq horizontal'' subbundle $H\subset \bT$. In the first approach a splitting of the sequence (3) gives a map $\varpi:J(E)\rightarrow E\otimes T^{*}X$. From the definition of the jet space, a section $s$ of $E$ induces a section $J(s)$ of $J(E)$ and $\nabla s= \varpi(J(s))$ is a covariant derivative on sections of $E$. While $C^{\infty}$ connections always exist {\it holomorphic} connections  need not and one of the  main results of Atiyah in \cite{kn:MFA2} is:
\begin{prop}
A bundle $E$ admits a holomorphic connection if and only if the Atiyah class
$\alpha(E)\in H^{1}({\rm End} E)$ is zero. 
\end{prop}

For another point of view on the Atiyah class we consider a $C^{\infty}$ connection, or covariant derivative,  $\nabla_{E}$ on $E$ which is compatible with the holomorphic structure in the sense that the $(0,1)$ component of $\nabla_{E}$ equals $\db_{E}$. (In the principal bundle approach this is the same as saying that the horizontal subspaces are complex subspaces of ${\cal T}$.) The connection has curvature $F_{\nabla}\in \Omega^{2}({\rm End} E)$ and since $\db_{E}^{2}=0$ the (0,2) part of $F_{\nabla}$ vanishes and
$  F_{\nabla}= F_{\nabla}^{1,1} + F_{\nabla}^{2,0} $.
The Bianchi identity implies that $\db_{{\rm End} E} F_{\nabla}^{1,1}=0$, so $F_{\nabla}^{1,1}$ defines a class in $H^{1}({\rm End} E \otimes T^{*}X)$ which is another representation of the Atiyah class.

This makes a connection with Chern-Weil theory.  Let $\psi$ be an invariant polynomial of degree $k$ on the Lie algebra of $G=GL(r,\bC)$. (In other words, $\psi$ is a polynomial function on matrices with $\psi(gmg^{-1})=\psi(m)$: these are just symmetric functions in the eigenvalues.) Then the combination of $\psi$ and wedge product defines a  map from the tensor product of $k$ copies of ${\rm End}\ E \otimes T^{*}X$ to $\Lambda^{k}T^{*}X$. If we have any class $\beta\in H^{1}({\rm End}\ E\otimes T^{*}X)$ we combine the map above with the product on cohomology to get  $\psi(\beta)\in H^{k}(\Lambda^{k} T^{*}X)$. Suppose, for simplicity, that $X$ is a compact K\"ahler manifold. Then we have a Hodge decomposition of the cohomology and
$$   H^{k}(\Lambda^{k}
T^{*}X)=H^{k,k}(X)\subset H^{2k}(X,\bC). $$ 
So $\psi(\beta)$ can be viewed as a class in the topological cohomology of $X$. Applying this to $\beta=\alpha(E)$  we get  the usual Chern-Weil construction for characteristic classes of $E$. This is clear from the representation of the Atiyah class by the curvature. We can choose a connection compatible with a Hermitian structure on the bundle. The curvature has type $(1,1)$ and writing out the recipe above, using the Dolbeault description where the product on cohomology is induced by wedge product on forms, gives exactly the Chern-Weil construction. (The whole discussion extends to  general
structure groups $G$ and that is  the context in which Atiyah worked.)

The curvature of a holomorphic connection has type $(2,0)$, so  
when $X$ is a Riemann surface a holomorphic connection is flat.
Using this,  Atiyah obtained alternative proofs of  some results of Weil from his 1938 paper \cite{kn:Weil} which began the study of higher rank vector bundles over Riemann surfaces. In the classical case of line bundles of degree $0$ a standard approach is to view sections as automorphic functions on the universal cover, transforming according to a multiplier. In other words, the line bundle is endowed with a flat connection with structure group $\bC^{*}$ and arises from a representation $\rho:\pi_{1}(X)\rightarrow \bC^{*}$. Weil extended this to higher rank bundles and flat bundles with structure group $GL(r,\bC)$. He showed that an indecomposable bundle $E$ over a compact Riemann surface admits a flat connection if and only if it has degree (i.e. first Chern class) equal to $0$. For a bundle $E$ over $X$ we have the Serre duality:
$$  H^{1}({\rm End} E\otimes T^{*}X)= H^{0}({\rm End E})^{*}. $$
So a flat connection exists if any only if the pairing between the Atiyah class and every holomorphic section of ${\rm End} E$ vanishes.
If the bundle is indecomposable (i.e. cannot be written as a non-trivial direct sum) then any holomorphic section of ${\rm End} E$ can be written as $\lambda 1_{E} + N$, where $N$ is nilpotent. (For otherwise the eigenspaces would give a decomposition.) By the discussion  of Chern-Weil theory above, the pairing between $1_{E}$ and $\alpha(E)$ is the degree of $E$,  so one has to see that the pairing of $\alpha(E)$ with the nilpotent sections vanishes. Suppose, for simplicity that $E$ has rank $2$ and there is a rank $1$ subbundle
$L= {\rm ker}\ N= {\rm Im}\ N$. Then we can choose a connection $\nabla$ on $E$ that preserves $L$. The pairing is given by
$$  \int_{X} {\rm Tr} (F_{\nabla}\ N), $$
and the integrand vanishes since $F_{\nabla}$ preserves $L$. 

\subsection{Bundles over curves}

In this subsection we discuss the papers \cite{kn:MFA1} and \cite{kn:MFA3} of Atiyah which are more detailed investigations of holomorphic bundles, mainly over algebraic curves (or compact Riemann surfaces). 

The first paper \cite{kn:MFA1} is focused on vector bundles $E\rightarrow X$ of rank $2$. Then the projective bundle $\bP(E)$ is a ruled surface so the study becomes part of the general study of algebraic surfaces. The starting point is the fact that any such bundle contains rank $1$ subbundles $L\subset E$.  For if $H$ is any positive line bundle over $E$ and $k$ is sufficiently large the tensor product  $E\otimes H^{k}$ has a non-zero holomorphic section. This section  might vanish at some points of $X$ but it is still true that the image of the section defines a rank $1$ subbundle of $E\otimes H^{k}$.  This is a special feature of  dimension $1$: in a local trivialisation about a point $x$ in $X$ the section is defined by a vector-valued function and we take the fibre of the subbundle over $x$ to be the line spanned by the first non-vanishing derivative.  This  fact means that we can write $E$ as an extension
\begin{equation}  0\rightarrow L\rightarrow E \rightarrow Q\rightarrow 0 \end{equation}
where $L$ and $Q$ are line bundles and by the theory discussed in the previous subsection the bundle $E$ is determined by $L, Q$ and an extension class in $H^{1}(L\otimes Q^{*})$. Multiplying the extension class by a non-zero scalar just corresponds to scaling the inclusion map $i:L\rightarrow E$  and does not change the rank 2 bundle. So if we leave out the trivial extension the data is $(L,Q,\xi)$ where $\xi\in \bP(H^{1}(L\otimes Q^{*})\setminus \{0\}$. Thus, using the theory of extensions, Atiyah was able to \lq\lq reduce'' questions about rank 2 bundles to questions about line bundles and cohomology. The essential difficulty is that the line subbundle is not unique, so any bundle $E$ has  a multitude of descriptions of this form.

To reduce this multitude, Atiyah considers subbundles of  {\it maximal} degree.
As we will see, this is an idea that goes a long way, so we will make a definition.
\begin{defn}
For a rank 2 bundle $E$ over a compact Riemann surface $X$  let $\nu(E)$ be the maximal degree of a rank 1 subbundle of $E$.
\end{defn}
(It is straightforward to show that this is well-defined.)

Let us illustrate the utility of this notion by classifying rank 2  bundles over the  Riemann sphere $\bP^{1}$. Recall that the line bundles over $\bP^{1}$ are just the powers $\cO(k)$ of the Hopf bundle $\cO(1)$ and that $H^{0}(\cO(k))=0 $  for $k< 0$ while $H^{1}(\cO(k))=0$  for $k> -2$. By taking the tensor product with a line bundle we may reduce to the case when $\nu(E)=0$, so $E$ is an extension 
\begin{equation}  0\rightarrow \cO\rightarrow E\rightarrow \cO(-k)\rightarrow 0 \end{equation}
 but $E$ has no subbundle of strictly positive degree.  Take the tensor product of (6) with $\cO(-1)$ and the long exact sequence in cohomology, which runs

$$ \dots H^{0}(E(-1))\rightarrow H^{0}(\cO(-k-1))\rightarrow H^{1}(\cO(-1))\dots $$
(Here, and later, we write $E(p)$ for $E\otimes \cO(p)$.)
If $k\leq - 1$  there is a section of $\cO(-k-1)$ which lifts to a section of $E(-1)$ since $H^{1}(\cO(-1))=0$. The image of this section is a line subbundle $L$ of $E$ and the section can be written as a section of $L(-1)$ composed with inclusion map, so the degree of $L$ is strictly positive which contradicts our hypothesis. Thus we see that $k>-1$. Now the extension class of (6) lies in $H^{1}(\cO(k))$ which vanishes for $k>-1$, so we conclude that the extension splits and $E=\cO\oplus \cO(-k)$. (The definition of $\nu$ implies that in fact $k\geq 0$.) In terms of ruled surfaces the result states that the $\bP^{1}$ bundles over $\bP^{1}$ are the {\it Hirzebruch surfaces} $\Sigma_{k}=\bP(\cO\oplus \cO(-k))$ for $k\geq 0$,  which Atiyah refers to  as the classical result \lq\lq that every rational normal ruled surface can be generated by a 1-1 correspondence between two rational normal curves lying in skew spaces''.

One of Atiyah's main results in \cite{kn:MFA1} determines when two extensions define isomorphic rank 2 bundles over a Riemann surface. Start with one description which we can take to be
\begin{equation}  0\rightarrow \cO\rightarrow E\rightarrow L_{1}\rightarrow 0, \end{equation}
for a line bundle $L_{1}$ of degree $d$ and an extension class in $ H^{1}(L_{1}^{*})$. Assume that the extension does not split, so $E$ is determined by a point $\eta$  in the projective space  $\bP(H^{1}(L_{1}^{*}))$ which we denote by $ \bP_{L_{1}}$. The Serre dual of $H^{1}(L_{1}^{*})$ is $H^{0}(L_{1}\otimes K_{X})$ so we have $\bP_{L_{1}}= \bP(H^{0}(L_{1}\otimes K_{X})^{*})$ and the linear system $L_{1}\otimes K_{X}$ defines a map $f: X\rightarrow \bP_{L_{1}}$. 

Now take a non-trivial line bundle $\Lambda$ of degree $0$ and consider the possible existence of an extension 
$$   0\rightarrow \Lambda \rightarrow E\rightarrow L_{1}\otimes \Lambda^{-1}\rightarrow 0, $$
(where the last term is fixed by considering the determinant of $E$). Taking a tensor product with $\Lambda^{-1}$, this is equivalent to considering extensions
$$   0\rightarrow \cO\rightarrow E\otimes \Lambda^{-1}\rightarrow L_{2}\rightarrow 0, $$
where $L_{2}= L_{1}\otimes \Lambda^{-2}$. The existence of such an extension is equivalent to a non-vanishing section of $E\otimes \Lambda^{-1}$ and from (7)  we have an exact sequence
$$ 0\rightarrow   H^{0}(E\otimes \Lambda^{-1})\rightarrow H^{0}(L_{1}\otimes \Lambda^{-1})\rightarrow H^{1}(\Lambda)\dots . $$
So we need a section of $L_{1}\otimes \Lambda^{-1}$ which maps to zero in the $(g-1)$-dimensional vector space $H^{1}(\Lambda)$. A section of $L_{1}\otimes \Lambda^{-1}$ defines a positive divisor say  $D=Q_{1}+\dots Q_{d}$ for points $Q_{i}\in X$. Now we have $L_{1}\otimes L_{2}= L_{1}\otimes L_{1}\otimes \Lambda^{-2}= (L_{1}\otimes \Lambda^{-1})^{2}$. So 
\begin{equation} \cO(2D)= L_{1}\otimes L_{2}. \end{equation}
Atiyah shows that this section maps to $0$ in $H^{1}(\Lambda)$ if and only if {\it the point $\eta\in \bP_{L_{1}}$ is in the linear span of $f(Q_{1}),\dots f(Q_{d})$}. So, starting with $L_{1}$ and $\eta\in \bP_{L_{1}}$, we consider the curve $f(X)\subset \bP_{L_{1}}$  and find all the configurations of $d$ points on this curve whose linear span contains $\eta$. For each such configuration we define another line bundle $L_{2}$ by (8) and this gives all the  descriptions of the rank 2 bundle by extensions with a sub-bundle of maximal degree.

Atiyah used these techniques in \cite{kn:MFA1} to classify all rank 2 bundles over curves of genus $1$ and $2$, extending the classical case of genus $0$ that we saw above. Postponing the discussion of genus $1$, we consider now  genus $2$. There are various cases, depending on the first Chern class of the bundle and the invariant $\nu$. By taking the tensor product with a line bundle it suffices to consider bundles $E$ with $c_{1}(E)$ equal to $0$ or $1$. The most interesting case, for reasons we will see more of later, is when $c_{1}(E)=1$ and $\nu(E)=0$. So, after tensoring by a line bundle, we have an extension (7) where the line bundle $L$ has degree $1$. Thus $L\otimes K_{X}$ has degree $3$; the Riemann-Roch formula shows that $H^{0}(L\otimes K_{X})$ has dimension $2$ and $\bP_{L}$ is a projective line. Since the dimension of this vector space is $2$ we have a canonical identification between the projective space and its dual so, in this special situation, we can regard the extension data as the  zero divisor $P+Q+R$ of a section of $L\otimes K_{X}$. Conversely, any divisor $D=P+Q+R$ defines a degree $1$ line bundle $L=\cO(D)\otimes K_{X}^{-1}$ and an extension class. So for each triple of points $P,Q,R$ we have a rank $2$ vector bundle $E(PQR)$ say. The curve $X$ is hyperelliptic, with an involution $\tau:X\rightarrow X$. Applying his criterion, Atiyah shows that a different triple $P'Q'R'$ defines a bundle $E(P'Q'R')$ projectively equivalent to $E(PQR)$ if and only if the new triple is given by applying $\tau$ to {\it two} of the points: for example  $P'=P, Q'=\tau(Q), R'=\tau(R)$.  In this way he constructs a $3$-dimensional moduli space $M_{\bP}$ parametrising these projective bundles.

 To see this moduli space more explicitly we start with the symmetric product
 $s^{3}(X)$ of triples $PQR$. This is a fibre bundle over the Jacobian $\Jac(X)$ with fibre $\bP^{1}$. There is a $4-1$ map from $s^{3}(X)$ to $M_{\bP}$ given by the identifications described above. This map takes the $\bP^{1}$ fibres in the symmetric product to the bundles constructed using a fixed line bundle but varying extension data. We have a $2-1$ map from $X$ to $\bP^{1}$ which induces an $8-1$ map from $s^{3}(X)$ to $s^{3}(\bP^{1})=\bP^{3}$. This map factors through the space $M_{\bP}$ so we have
$$  s^{3}(X)\rightarrow M_{\bP}\stackrel{\pi}\rightarrow \bP^{3}, $$
where $\pi$ is $2-1$. We have $6$ branch points $\lambda_{i}\in \bP^{1}$ of
the double covering $X\rightarrow \bP^{1}$. Each of these points $\lambda_{i}$ defines a plane $\Pi_{i}$ in $\bP^{3}=s^{3}(\bP^{1})$. The 3-fold $M_{\bP}$ is the double cover of $\bP^{3}$ branched over these six planes $\Pi_{i}$.

\

We now turn to Atiyah's results on bundles over elliptic curves (genus $1$). The rank 2 case was covered in \cite{kn:MFA1} and the later paper \cite{kn:MFA3} gave a complete classification for all ranks. Slightly before, Grothendieck classified bundles of all ranks over $\bP^{1}$: they are all direct sums of line bundles and the proof follows the same lines as the rank $2$ case discussed above, using  induction on the rank.

Let $X$ be a compact Riemann surface of genus $1$ and fix a line bundle $\lambda$ of degree $1$ over $X$.
For integers $r,d$ with $r>0$ let $\cE(r,d)$ be the set of isomorphism classes of indecomposable bundle of rank $r$ and degree $d$. The most straightforward case is when $r,d$ are coprime so we begin with that. Atiyah proved
\begin{prop}
There is a unique way to define indecomposable bundles $E_{r,d}$ (up to isomorphism) for coprime rank $r$ and degree $d$ such that
\begin{itemize}
\item $E_{1,0}=\ubC$
\item $E_{r,d+r}= \lambda\otimes E_{r,d}$
\item If $0<d<r$ there is an exact sequence
$$  0\rightarrow \ubC^{d}\rightarrow E_{r,d}\rightarrow \rightarrow E_{r-d,d}\rightarrow 0. $$
\end{itemize}

\end{prop}

Atiyah's construction of the $E_{r,d}$ follows the Euclidean algorithm for the pair $(r,d)$. To illustrate this, consider the case of $E_{5,3}$. Following the third bullet,  we want to build this as an extension
\begin{equation} 0\rightarrow \ubC^{3}\rightarrow E_{5,3}\rightarrow E_{2,3}\rightarrow 0, \end{equation} so we first have to construct $E_{2,3}$. Following the second bullet, we have $E_{2,3}=\lambda \otimes E_{2,1}$ so we have to construct $E_{2,1}$. The three bullets together state that this  should be an extension
\begin{equation}  0\rightarrow \ubC\rightarrow E_{2,1}\rightarrow \lambda \rightarrow
0. \end{equation}
These extensions (10) are classified by $H^{1}(\lambda^{*})$ which is 1-dimensional, so there is indeed a unique indecomposable $E_{2,1}$. Standard arguments show that, with $E_{2,3}=E_{2,1}\otimes \lambda$,  the cohomology $H^{1}(E_{2,3}^{*})$ is $3$-dimensional. So an isomorphism $h:H^{1}(E_{2,3}^{*})\rightarrow \bC^{3}$ builds an extension (9) and the middle term $E_{5,3}$ does not depend (up to isomorphism) on the choice of $h$ (for any two choices differ by an automorphism of $\bC^{3}$).

Now the basic fact is that this construction essentially gives all indecomposable bundles with $(r,d)$ coprime: as $L$ runs over the line bundles of degree $0$ the tensor products $E_{r,d}\otimes L$ run  once over all of $\cE(r,d)$, so $\cE(r,d)$ is identified with the Jacobian $\Jac(X)$.

The statements for $(r,d)=h>1$ are similar except that we start with a more complicated initial bundle $F_{h}\in \cE(h,0)$. These are the unique indecomposble bundles with degree $0$, rank $h$ and having  non-trivial sections and they can be built inductively as
$$   0\rightarrow \bC\rightarrow F_{h}\rightarrow F_{h-1}\rightarrow 0. $$
In the end, each $\cE(r,d)$ is a copy of the Jacobian. 

\subsection{Double points in dimension 2 and 3}

The paper \cite{kn:MFA4} contains some of the most significant results from Atiyah's early work. It lies a little outside the main themes that we have discussed but vector bundles  have some  bearing on the material.

  Consider a hypersurface $X$ in $\bC^{n+1}$ defined by an equation $f(z)=0$ with one singular point at the origin, so   the analytic function $f$ and its first derivatives vanish at $z=0$. There are two operations we can perform to replace $X$ by a complex manifold:
\begin{itemize}
\item Resolution:  blowing up the origin in $\bC^{n+1}$ sufficiently many times with suitable centres, the proper transform of $X$ is a  complex manifold $\tX$ with a map $\pi:\tX\rightarrow X$ and an exceptional divisor $E\subset \tX$ which is collapsed to the origin. The map gives a holomorphic isomorphism from  $\tX\setminus E$ to $X\setminus\{0\}$. . 
\item Smoothing: for small non-zero $\delta\in \bC $ we get a smooth hypersurface $X_{\delta}$ defined by the equation $f(z)=\delta$.
\end{itemize}

The spaces $\tX$ and $ X_{\delta}$ are   usually  completely different as smooth manifolds.  The singularity is called an ordinary double point if the Hessian $\frac{\partial^{2}f}{\partial z_{i}\partial z_{j}}$of $f$ at the origin is nonsingular. Atiyah's first result in \cite{kn:MFA4} is that:
\begin{prop} For ordinary double points of complex surfaces the manifolds $\tX, X_{\delta} $ are diffeomorphic.\end{prop}
 More generally if we have any surface $S$ with ordinary double point singularities which arises as a limit of smooth surfaces $S_{\delta}$ then the resolution $\tS$ is diffeomorphic to the smoothings $S_{\delta}$.

If the Hessian of $f$ is nonsingular then, as Atiyah proved, there is a local holomorphic co-ordinate transformation taking $f$ to the standard quadratic form $\sum z_{i}^{2}$. (This is a holomorphic version of the Morse lemma for smooth real-valued  functions.) In this way one can reduce to the case when $X$ is the singular quadric surface $z_{1}^{2}+z_{2}^{2}+z_{3}^{2}=0$. To see how this is related to vector bundles, consider extensions over $\bP^{1}$:
  \begin{equation}  0\rightarrow \cO(-1)\rightarrow E \rightarrow \cO(1)\rightarrow 0, \end{equation}
classified by $H^{1}(\cO(-2))=\bC$. For each point $\tau\in H^{1}(\cO(-2))$ we can construct a bundle $E_{\tau}$ and these fit into a smooth family as $\tau$ varies and  it is clear that the projective bundles $\bP(E_{\tau})$ are all diffeomorphic. When $\tau$ is non-zero the bundle $E_{\tau}$ is holomorphically trivial: the sequence (11) is just that arising from the inclusion of the tautological bundle $\cO(-1)$ in the trivial bundle $\ubC^{2}$. Thus for $\tau\neq0$ the ruled surface $\bP(E_{\tau})$ is the product $\bP^{1}\times \bP^{1}$, the Hirzebruch surface $\Sigma_{0}$. When $\tau$ is zero
we have $E_{\tau}=\cO(1)\oplus \cO(-1)$ which has the same projectivisation as $\cO\oplus \cO(2)$ so  $\bP(E_{0})$ is the Hirzebruch surface $\Sigma_{2}$. We see that $\Sigma_{2}$  is diffeomorphic to $\Sigma_{0} = \bP^{1}\times\bP^{1}$ but with a different complex structure. In general, the ruled surfaces $\Sigma_{k}$ for even $k$ are all diffeomorphic (as Atiyah mentions in \cite{kn:MFA1}).

For $\delta\neq 0$ the projective completion in $\bP^{3}$ of the  affine surface $\sum z_{i}^{2}=\delta$ is a non-singular projective quadric surface $S_{\delta}$, hence equivalent to $\bP^{1}\times \bP^{1}$. The projective completion $S_{0}$ of the singular affine surface
$\sum z_{i}^{2}=0$ is a cone over the conic $Q$ in the hyperplane at infinity defined by the same equation. Blowing up the origin,  we get a resolution $\tS_{0}\rightarrow S_{0}$ with exceptional divisor $E$ which is a copy of $Q$ with self-intersection $E.E=-2$.   The proper transforms of the lines in $\tS_{0}$ through the origin give a ruling of $\tS_{0}$ which displays   $\tS_{0}$ as the Hirzebruch surface $\Sigma_{2}$. So the preceding discussion shows that $\tS_{0}$ is diffeomorphic to $S_{\delta}$ and it is not hard to see that this implies  Proposition 4.  

 In fact there are many ways of understanding this.  If we take $\delta$ to be a small positive real number then the complex quadric $X_{\delta}$ in $\bC^{n+1}$ defined by the equation $\sum z_{i}^{2}=\delta$ contains an $n$-dimensional sphere as the set of real solutions. It is easy to show that, as a $C^{\infty}$ manifold, $X_{\delta}$ fibres over $S^{n}$ with fibres $\bR^{n}$. On the other hand the resolution $\tX$ contains the exceptional divisor which can be identified with a complex quadric $Q\subset \bC\bP^{n}$ and  $\tX$ is a holomorphic fibration  over $Q$ with fibre $\bC$. For $n>2$ the spaces $X_{\delta}$ and $\tX$ do not even have the same homotopy type: the special feature of the surface case, when $n=2$, is that the conic $Q\subset \bP^{2}$ is diffeomorphic to the $2$-sphere.

Something special also happens when $n=3$ and this is the  most original and significant aspect of Atiyah's paper: the discovery of \lq\lq small resolutions''.
We can explain this, in a way which will fit in with material later in this article, by considering quotient spaces. Consider first the action of $\bC^{*}$ on $\bC^{2}$ given by 
$$  \lambda\mapsto \left(\begin{array}{cc} \lambda &0\\0&\lambda^{-1}\end{array}\right). $$
There is an invariant polynomial $xy$ and for $t\neq 0$ the set $xy=t$ is one orbit of the action.  But in the exceptional case $t=0$ there are 3 distinct orbits: the origin and the two axes minus the origin. The quotient space $\bC^{2}/\bC^{*}$ is a non-Hausdorff space since these orbits cannot be separated. To get around this we make a choice of which orbit to include and then we get a quotient space $\bC$. 

Now extend this to consider the action of $\bC^{*}$ on $\bC^{4}$ with co-ordinates $(x_{1}, x_{2}, y_{1}, y_{2})$ with $\lambda$ acting as multiplication by $\lambda$ on the $x_{i}$ and $\lambda^{-1}$ on the $y_{j}$. There are 4 invariant polynomials $A_{ij}= x_{i} y_{j}$ and these satisfy an identity
\begin{equation}  A_{11} A_{22}- A_{12} A_{21}=0. \end{equation}
In more invariant terms, if we write $\bC^{4}=V_{1}\otimes V_{2}$ for $2$-dimensional spaces $V_{i}$ then the polynomials $A_{ij}$ define the tensor product map  $\alpha: V_{1}\times V_{2}\rightarrow V_{1}\otimes V_{2}$ and the image is the variety of matrices with determinant zero. 

The left hand side of (12) is a nondegenerate quadratic form so the image of $\alpha$ can be identified with our standard singular quadric $X\subset \bC^{4}$ with one ordinary double point. 
 The fibre of $\alpha$ over any nonzero point is a single $\bC^{*}$ orbit but the quotient $\bC^{4}/\bC^{*}$ is not Hausdorff for the same reason as before. So, as before,  we consider three different subsets $U_{0},U_{+}, U_{-}$ of $\bC^{4}$:
\begin{itemize}
\item in $U_{+}$ we include all points $(x,y)$ with $x\neq 0$;
\item in $U_{-}$ we include all points $(x,y)$ with $y\neq 0$;
\item in $U_{0}$ we include the origin but no other point $(x,y)$ with $x$ or $y$ equal to $0$.
\end{itemize}
Then $U_{0}/\bC^{*}$ is identified by $\alpha$ with the singular quadric $X$ while $Y_{\pm}= U_{\pm}/\bC^{*}$ are complex manifolds. The projectivization map from $V_{1}\setminus \{0\}$ to $\bP(V_{1})$  induces a  map from  $Y_{+}$ to $\bP^{1}$ with fibre $\bC^{2}$. More precisely, $Y_{+}$ is the total space of the bundle $\cO(-1)\oplus \cO(-1)$ over $\bP^{1}$. The map $\alpha$ induces a map from $Y_{+}$ to the singular quadric $X$ with fibre $\bP^{1}$ over the origin. This is a {\it small} resolution because this fibre has codimension $2$ rather than $1$.  Of course the same applies to $Y_{-}$,  but this gives a {\it different} small resolution. In the usual resolution $\tX$, obtained by blowing up the origin,  the exceptional fibre is a copy of  a nonsingular quadric surface $Q\subset \bC\bP^{3}$ which is the product $\bP^{1}\times\bP^{1}$. Then the map $\tX
\rightarrow X$ factors through either $Y_{+}$ or $Y_{-}$ by collapsing one of the $\bP^{1}$ factors in $Q$ and  we have a diagram:

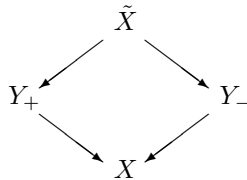
\begin{picture}(1000,320)(-300,0)
\put(500,300){$\tX$}
\put(500,20){$X$}
\put(300,160){$Y_{+}$}
\put(700,160){$Y_{-}$}
\put(680,140){\vector(-4,-3){120}}
\put(560,280){\vector(4,-3){120}}
\put(360,140){\vector(4,-3){120}}
\put(480,280){\vector(-4,-3){120}}
\end{picture}

To relate this to the two dimensional case we take the composition of the map $Y_{+}\rightarrow \bC^{4}$ with a generic linear projection from $\bC^{4}$ to $\bC$ to get $p: Y_{+}\rightarrow \bC$. For nonzero $t$  in $\bC$ the fibre $p^{-1}(t)$ is a  smooth quadric surface in $\bC^{3}$ while $p^{-1}(0)$ is the resolution of the singular quadric surface. One checks that $p$ is a submersion which shows again that the two are diffeomorphic, and this is the way that Atiyah proves Proposition 4.

Finally, suppose that $Z$ is a 3-fold and $q:Z\rightarrow B$ is a holomorphic map to a Riemann surface such that the fibres have at worst ordinary double point singularities. In other words the triple $(Z,B,q)$ is a family of surfaces over $B$. Atiyah considers a double cover  $\varpi:\hat{B}\rightarrow B$ branched over the critical values of $q$ (and possibly other points of $B$).  Then there is a pulled back family $\hat{Z}\rightarrow \hat{B}$ but each critical point of $q$ produces an ordinary double point singularity in $\hat{Z}$. Performing small resolutions of these singularities Atiyah gets a new family $q^{*}: Z^{*}\rightarrow \hat{B}$ which is a submersion with all fibres smooth. In other words, he resolves the singularities of the fibres of $q$ \lq\lq in the family'', after going to the double cover. (The cover cannot be avoided since the original family has non-trivial monodromy around the critical values in $B$ but the monodromy has order $2$ so is removed in the double cover.)

\subsection{Some later developments}

In Part II we jump forward almost two decades so here we will mention some developments in the intervening period, roughly 1958-1977, which are relevant to both Atiyah's early and later work. 

   The theory of moduli spaces of bundles over Riemann surfaces, which had appeared in an informal way in Atiyah's early papers, was put on solid foundations in the early 1960's by Mumford \cite{kn:Mumford}, as an application of his Geometric Invariant Theory. He introduced the notion of \lq\lq stability'' of a vector bundle over a compact Riemann surface $X$ and showed that each $(r,d)$ there is a moduli space $M(r,d)$ of stable holomorphic or rank $r$ and degree $d$ over $X$. These map to  the Jacobian of $X$ by the determinant of the bundle and the fibre, $N(r,d)$ say, has complex dimension $(r^{2}-1)(g-1)$.
These are quasiprojective varieties and when $r$ and $d$ are coprime they are projective. In the case of rank $2$ bundles the stability condition is just that $\nu(E)< d/2$ where $\nu$ is the invariant we discussed in 1.2 above. The geometry and topology of these moduli spaces was studied by Newstead (a former student of Atiyah), Narasimhan, Seshadri,  Ramanan and others. The simplest case is when $X$ has genus $2$ with $r=2,d=1$. It was shown by Newstead \cite{kn:Newstead2} and Narasimhan and Ramanan \cite{kn:MSNR} that $N=N(2,1)$ is the intersection of two quadrics in $\bP^{5}$. Up to a covering this is the same as  Atiyah's moduli space of projective bundles, discussed in subsection 1.2 above. If $L$ is a line bundle over $X$ with $L^{2}=\ubC$ then  rank $2$ bundles $E$ and $E\otimes L$ have the same determinant. In genus $2$ there are $2^{4}$ such line bundles $L$ and we get an action of $(\bZ/2)^{4}$ on $N$ with quotient $M_{\bP}$. Write $N$ in standard form as the intersection
of the quadrics $\sum w_{i}^{2}=0$ and $\sum \lambda_{i} w_{i}^{2}=0$ in $\bC\bP^{5}$. The group $G=(\bZ/2)^{6}$ acts on $\bC^{6}$ by multiplication of the co-ordinates $w_{i}$ by $\pm 1$. Let $G_{0}\subset G$ be the index $2$ subgroup which changes only an even number of signs. Then $-1$ is in $G_{0}$ and acts trivially on $\bP^{5}$ and $G_{1}=G_{0}/(\pm 1)$ is isomorphic to $(\bZ/2)^{4}$. This gives the $(\bZ/2)^{4}$ action on $N$ and $M_{\bP}= N/G_{1}$. To relate this to Atiyah's description, the map $Z_{i}=w_{i}^{2}$ induces an equivalence between $N/G$ and the 3-plane $\bP^{3}=\{\sum Z_{i}=0, \sum \lambda_{i} Z_{i}=0\}$ in $\bC\bP^{5}$. The induced map from $M_{\bP}= N/G_{1}$ to $\bP^{3}$ has degree $2$ (since $\pm 1$ acts trivially on $N$) and is a double cover branched over the six planes $\{ Z_{i}=0\}$ in $\bP^{3}$. 
We saw that the extension construction gives a map from a $\bP^{1}$ bundle over the Jacobian $\Jac(X)$ to $M_{\bP}$. This lifts to $N$ if we take  the cover of $\Jac(X)$ by itself given by the doubling map. So for each point in $\Jac(X)$ we get a map from $\bP^{1}$ to $N$.  Newstead and Narasimhan and Ramanan showed that the images are lines and   all  lines in $N$ arise in this way.  This fits into a famous classical picture involving three spaces (see \cite{kn:GH} Chapter 6).
\begin{itemize} \item  Genus $2$ curves $X$.
\item Principally polarised abelian surfaces $J$.
\item Intersections $N$ of two quadrics in $\bP^{5}$.
\end{itemize} 

To get from $X$ to $J$ we take the Jacobian $\Jac(X)$. To get from $J$ to $X$ we take the theta divisor (assuming that $J$ is not a product). To get from $X$ to $N$ we take the moduli space of rank $2$ bundles $N(2,1)$. To get from $N$ to $X$ we take the six points $\lambda_{i}$ in $\bP^{1}$ corresponding to the singular quadrics in the pencil through $N$ and the double cover of $\bP^{1}$ branched at these points. Finally, to get from $N$ to $J$ we take the set of lines in $N$.

In subsection 1.1 we mentioned the 1938 result of Weil which states that an indecomposable  rank $r$ vector bundle of degree $0$ over a curve arises from a representation of the fundamental group in $GL(r,\bC)$.  The drawback of this is that the representation is far from unique. For example in the case of rank 1 bundles the space of representations has dimension $2g$ while the space of degree $0$ line bundle---the Jacobian---has dimension $g$. In 1965 Narasimhan and Seshadri obtained a much sharper result. First, the degree restriction can be removed by considering projective representations. More important they showed that considering irreducible unitary representations selects out exactly the stable bundles of Mumford and the representation associated to a bundle is unique. More precisely, let $x_{0}\in X$ be a base point and $\alpha_{i},\beta_{i}$ be standard generators for the fundamental group $\pi_{1}(X,x_{0})$ with relation 
$$  \prod_{i=1}^{g} [\alpha_{i},\beta_{i}] =1. $$
Let $\tilde{\pi}_{r}$ be the central extension of $\pi_{1}(X,x_{0})$ defined by adjoining an element $z$ such that $z^{r}=1$ and 
$\prod_{i=1}^{g} [\alpha_{i},\beta_{i}] =z$. 
Then Narasimhan and Seshadri proved in \cite{kn:MSNCS} that there is a 1-1 correspondence between equivalence classes of  stable rank $r$  bundles $E$ with $\Lambda^{r}E= d\ \cO(x_{0})$ and
conjugacy classes of irreducible unitary representation of $\tilde{\pi}_{r}$ in $U(r)$ which map $z$ to ${\rm exp}(2\pi i d/r)$.   It is hard to overstate the significance of this result for later developments, some of which we will discuss in Part II. Here we will just see how it fits in with Atiyah's classification of bundles over elliptic curves. In that case stable bundles only occur when $r,d$ are coprime, so we assume that. The representations are given by $A,B\in U(r)$ satisfying
  $$  ABA^{-1}B^{-1}= \zeta 1$$
  where $\zeta= {\rm exp}(2\pi id/r) \in \bC$ is a primitive rth. root of unity. To construct such a representation, let $A$ be the map which permutes the standard basis vectors cyclically
$$  A e_{1}=e_{2}, \dots , Ae_{r}= e_{1}, $$
and let $B$ be diagonal in this basis with $B e_{k}=\zeta^{k} e_{k}$.  This is the unique irreducible representation up to conjugacy  and  corresponds to Atiyah's bundle $E_{r,d}$.

 Meanwhile, there were developments in the study of bundles over higher dimensional varieties. Maruyama constructed moduli spaces of stable bundles over projective surfaces, with an extension of Mumford's definition of stability. Atiyah's former student Schwarzenberger studied bundles over the projective plane and other surfaces \cite{kn:RLES}. He introduced the notion of \lq\lq jumping lines''. A stable bundle with first Chern class $0$ over the plane is trivial on generic lines but not on all. There is a  curve in the dual plane parametrising lines such the restriction is non-trivial--typically $\cO(1)\oplus \cO(-1)$. This jumping curve is an interesting invariant of the bundle. The   degree of the curve is equal to the second Chern class of the bundle. For bundles of rank 2 and $c_{2}=2$ we get a conic and the bundle can be reconstructed from the conic. One takes a suitable line bundle $L\rightarrow Q$ over the quadric surface $Q$ which is a double cover of $\bP^{2}$ branched over the conic. The direct image of this line bundle is a rank 2 bundle over $\bP^{2}$.

The ramifications of Atiyah's results on double points are so extensive that they would require much more space then we have and much more knowledge on the part of the author to treat properly. In one direction they lead on to Brieskorn's theory of simple singularities associated to the Lie algebras of types A,D,E (which we will see a little of in 2.4.3 below). In another direction, small resolutions are  crucial in the classification theory of higher dimensional varieties and in Mirror symmetry for Calabi-Yau threefolds.

\section{ Gauge theory }

In the commentary on Volume 5 of his collected works \cite{kn:MFA5}, Atiyah wrote:
{\it From 1977 onwards my interests moved in the direction of gauge theories and the interaction between geometry and physics....the stimulus came from two sources. On the one hand Singer told me about the Yang-Mills equations which through the influence of Yang were just starting to percolate into mathematical circles... The second stimulus came from the presence in Oxford of Roger Penrose and his group.}

\subsection{The Yang-Mills equations and self-duality}

{\it Gauge theory} refers to the study of a connection $A$ on a principal bundle $P$ with structure group $G$ over a manifold $M$. The curvature $F(A)$ is a $2$-form with values in the bundle ${\rm ad} \ P$ associated to the adjoint action of $G$ on its Lie algebra. If $M$ has a Riemannian or pseudo-Riemannian metric $g$ and if we have an invariant quadratic form on the Lie algebra we can write down the Lagrangian
\begin{equation}  \cL =\int_{M} \vert F(A)\vert^{2} d\mu, \end{equation} and associated Euler-Lagrange equations
\begin{equation}   d_{A}^{*} F(A)=0.  \end{equation}
This is a second order partial differential equation for the connection $A$. (Of course the Euler-Lagrange equations are generated by the formal expression (13)---there is no assumption that the integral is well-defined.) If the structure group $G$ is $S^{1}$ the curvature is an ordinary 2-form. If $M$ is Lorentz space-time $\bR^{3,1}$ this is the setting for classical electro-magnetism, with $F$ the electromagnetic field. Writing this in terms of a space-time splitting as a pair $(E,B)$ of electric and magnetic fields the Lagrangian is
$$  \int_{\bR^{3,1}} \vert E\vert^{2}-\vert B\vert^{2}, $$
and the equation $d^{*}F=0$ along with $dF=0$ (which is true for any curvature $2$-form) are the source-free Maxwell equations for $(E,B)$. In all that follows we will be concerned with the case when $M$ has a Riemannian metric $g$ and we usually suppose that $G$ is compact and the form on the Lie algebra is positive definite. Then the integral (13) is the square of the standard $L^{2}$ norm of the curvature, as the notation suggests. The passage from Lorentzian to Euclidean signature becomes relevant to physics in  Quantum Field Theory, which we will not attempt to say anything about here.

There are several special features of Yang-Mills theory when the dimension of the base manifold $M$ is $4$. The integral (13) and the Yang-Mills equation (14) are conformally invariant and there are first-order \lq\lq instanton'' equations whose solutions satisfy the second order Yang-Mills equations, in much the same way as  holomorphic functions on $\bC$ are harmonic. We assume that the Riemannian $4$-manifold $M$ is oriented; then there is a $*$-operator
$*:\Lambda^{2}\rightarrow \Lambda^{2}$ with $*^{2}=1$ and the $2$-forms decompose into the $\pm 1$ eigenspaces $\Lambda^{2}=\Lambda^{2}_{+}\oplus \Lambda^{2}_{-}$ of {\it self-dual} and {\it anti self-dual} forms. Then we can write
$$  F(A)=F^{+}(A)+F^{-}(A) $$
and the \lq\lq instantons''  have $F^{+}=0$ (self-dual connections) or $F^{-}=0$ (anti self-dual connections). If $M=\bR^{4}$ with a \lq\lq space-time'' splitting $\bR^{4}=\bR^{3}\times \bR$  we can write the curvature as a pair $(E,B)$ and the instanton condition is $E=\pm B$. 

Another special feature of dimension $4$ is that there are topological invariants---characteristic numbers---which are also given by integrals of quadratic expressions in the curvature. Let $b$ be an invariant quadratic form on the Lie algebra. In Chern-Weil theory, for any connection $A$ on a $G$ bundle $P\rightarrow M$, one considers a $4$-form $b(F(A))$ on $M$. One shows that this is a closed form and that its de Rham cohomology class depends only on the bundle $P$ not on the choice of connection. It defines a characteristic class $\kappa_{b}$ of $P$ in $H^{4}(M;\bR)$.
When $M$ is a compact oriented $4$-manifold we get a characteristic number
$$  \kappa_{b}(P)= \int_{M} b(F) . $$
Take $b$ to be the positive definite form used to define the Lagrangian. Then from the definition of the $*$-operator we have
$$  \kappa_{b}(P) =\int_{M} \vert F^{+}\vert^{2}  - \vert F^{-}\vert^{2}\  d\mu, $$
while
$$   \cL(A)= \int_{M} \vert F^{+}\vert^{2} + \vert F^{-}\vert^{2}  \ d\mu. $$

It follows that an instanton connection on $P$ is an absolute minimum of $\cL$: it is a $\pm$ self-dual connection depending on the sign of $\kappa_{b}(P)$. For example consider the group $G=SU(r)$ with the standard positive definite form $b(\xi)=-{\rm Tr}(\xi^{2})$ on its Lie algebra. Then Chern-Weil theory shows that
$$  \kappa_{b}(P)= - 8\pi^{2} \langle c_{2}(P), [M]\rangle, $$
where $c_{2}$ is the second Chern class. So on a bundle with $c_{2}<0$ we can look for self-dual connections and on a bundle with $c_{2}>0$ we look for ant-self-dual connections. If $c_{2}=0$ the formulae show that any instanton must be flat, with $F=0$.

Atiyah's first contribution to  Yang-Mills theory (with Hitchin and Singer) was an application of the Index Theorem to find the dimensions of moduli spaces of instantons. This was outlined  in a note \cite{kn:AHS1} with full details and applications in \cite{kn:AHS2}. For physics the main case of interest was when the base manifold is Euclidean $\bR^{4}$ which has the $4$-sphere $S^{4}$ as conformal compactification. Finite action solutions of the Yang-Mills equations over $\bR^{4}$ are equivalent to smooth solutions on $S^{4}$ (as was proved a few years later by Uhlenbeck). Changing orientation switches self-dual and anti self-dual. There is a better fit with standard conventions if one chooses to discuss anti-self-dual (ASD)  connections so we will follow that convention here (different to the convention in the papers of Atiyah, Hitchin, Singer). 
The first interesting case is the structure group $G=SU(2)$ and a bundle $P\rightarrow S^{4}$ with $c_{2}(P)>0$. A construction of t'Hooft (which we will discuss further below) gave families of solutions depending on approximately $5k$ parameters. Atiyah, Hitchin and Singer showed that the general solution depends on $8k-3$ parameters. More precisely, the moduli space $\cM_{k}$ of solutions is a manifold of dimension $8k-3$.

We will now review  Atiyah, Hitchin and Singer's analysis of the deformation problem in general. This involves variants of the ideas and technology used to study deformations of complex manifolds and holomorphic bundles.  Let $A$ be an ASD connection on a $G$-bundle $P\rightarrow M$ with $M$ a compact oriented $4$-manifold and $G$ compact. Recall that \lq\lq the difference of two connections is a tensor'': the general connection on $P$ can be written as $A+a$ where $a\in \Omega^{1}(\rad P)$:  the $1$-forms with values in the bundle $\rad P$. We have a coupled exterior derivative
$$ d_{A}: \Omega^{p}(\rad P)\rightarrow \Omega^{p+1}(\rad P)$$
and the formula for the curvature is
\begin{equation}  F(A+a)= F(A)+ d_{A} a + \frac{1}{2} [a,a], \end{equation}
where the last expression combines the bracket on the fibres of $\rad P$ with the wedge product on forms. Thus
\begin{equation} F^{+}(A+a)= d_{A}^{+} a + \frac{1}{2} [a,a]^{+}, \end{equation}

where the $+$ superscripts denote projection to the self-dual forms. The group $\cG$ of automorphisms of the bundle $P$ covering the identity on $M$ acts on the connections and connections in the same orbit are geometrically equivalent. By definition the moduli space ${\cal M}$ is the quotient of the set of solutions to the ASD equation by ${\cal G}$. We can write
\begin{equation}   g(A)= A- (d_{A} g )g^{-1} . \end{equation}
(Strictly,  the notation  here assumes that $G$ is a matrix group so that $g$ can be regarded as a section of an associated vector bundle.)  

For small $a$ we have, schematically, $F^{+}(A+a)= d^{+}_{A}a+ O(a^{2})$. If $g$ is close to the identity we can write $g=\exp(-u)$ for $u\in \Omega^{0}(\rad P)$ and, again schematically,
$g(A)= A+d_{A} u+ O(u^{2})$.  The linearised version of the deformation problem is expressed through a  complex
\begin{equation}  \Omega^{0}(\rad P)\stackrel{d_{A}}{\rightarrow} \Omega^{1}(\rad P) \stackrel{d^{+}_{A}}{\rightarrow} \Omega^{2}_{+}(\rad P), \end{equation}
with cohomology groups $H^{0}_{A}, H^{1}_{A}, H^{2}_{A}$. The cohomology $H^{1}_{A}$ consists of solutions of the linearised ASD equation modulo those that arise from gauge transformations of $A$,  so it is a natural candidate to be the tangent for the moduli space $\cM$ at the equivalence class $[A]$.  One of the main results is that this is true provided  that the connection $A$ is \lq\lq irreducible'' and the second cohomology $H^{2}_{A}$ vanishes. Here \lq\lq irreducible'' means that the bundle $P$ with connection $A$ does not reduce to any proper subgroup of $G$. (In fact one needs the weaker condition that the bundle does not reduce to a subgroup whose centralizer in $G$ is larger than the centre of $G$.) This implies that $H^{0}_{A}$ vanishes, because a non-zero covariant constant section of $\rad P$ induces a reduction of the bundle.

Index theory computes  the Euler characteristic $h_{0}-h_{1}+h_{2}$ of the complex (18). For this we consider, as standard in elliptic and Hodge theory, the operator
$$ D_{A}=d^{*}_{A}+ d^{+}_{A}: \Omega^{1}(\rad P)\rightarrow \Omega^{0}(\rad P)\oplus \Omega^{2}_{+}(\rad P), $$
whose kernel is isomorphic to $H^{1}_{A}$ and cokernel to $H^{0}_{A}\oplus H^{2}_{A}$, so the index of $D_{A}$ is minus the Euler characteristic. Suppose for a moment that $M$ is a spin manifold, so we have bundles $V^{\pm}\rightarrow M$ of positive and negative spinors. These are bundles with fibre $\bC^{2}$ and structure group $SU(2)$. We have
$  V^{+}\otimes V^{+}= \Lambda^{2} V^{+}\oplus S^{2} V^{+}$. The first summand has a fixed isomorphism with $\bC$ and the second with $\Lambda^{2}_{+}\otimes \bC$. On the other hand $V^{+}\otimes V^{-}$ is isomorphic to the complexified cotangent bundle of $M$. So the complexification of $D_{A}$ is an operator
$$  D_{A}^{\bC}: \Gamma( V^{-}\otimes E)\rightarrow \Gamma(V^{+}\otimes E), $$
where $E$ is the complex vector bundle $V^{+}\otimes \rad P$. This operator is the Dirac operator (from negative to positive spinors) coupled to the vector bundle $E$ so Atiyah, Hitchin and Singer conclude that
\begin{equation} {\rm ind} (D_{A}) = - \int_{M} {\rm ch}(V^{+}\otimes \rad P) \hat{A}(M) , \end{equation}

 (see the article of Dan Freed in this volume). From general considerations, the formula is still valid  if $M$ has no  spin structure.

In the case when $M=S^{4}$, Atiyah Hitchin and Singer used this index formula, to find exactly which bundles admit ASD connections. In this case we can assume that $G$ is  compact, simply connected and simple. The $G$-bundles over $S^{4}$ are classified by $\pi_{3}(G)$ and by Lie theory there is an $SU(2)$ subgroup $K\subset G$, unique up to conjugacy,  such that the inclusion defines an isomorphism $\bZ=\pi_{3}(SU(2))\rightarrow \pi_{3}(G)$. So the $G$-bundles are determined by an integer $k$ and this is compatible with reduction of structure groups. The existence result is:
\begin{thm} There exist  $G$ irreducible ASD connections $G$ connections over $S^{4}$ if and only if
\begin{itemize} \item for $G=SU(r)$, $k\geq r/2$;
\item for $G=Sp(r)$, $k\geq r$;
\item for $G=Spin(r)$ and $r\geq 7$, $k\geq r/4$;
 \item for $G=G_{2}$, $k\geq 2$;
 \item for $G=F_{4}, E_{6}, E_{7}, E_{8}$, $k\geq 3$.
 \end{itemize}
 \end{thm}
 The proof is, roughly speaking, to start with the known $SU(2)$ solutions and to compare moduli space dimensions to see when they can be deformed to larger structure groups. This involves  a more delicate analysis of the deformation theory, when $H^{0}_{A}$ is nonzero. 

\subsection{Self-dual manifolds and twistor spaces}

Another theme in the paper \cite{kn:AHS2} of Atiyah, Hitchin and Singer was a foundational study of self-dual 4-manifolds and their twistor spaces, in the Riemannian case. (Twistors were introduced by Penrose, initially over Lorentzian space-times.) 

Let $(M,g)$ be an oriented Riemannian manifold of even dimension $2m$. For each point $p\in M$ we can consider complex structures on the real vector space $TM_{p}$ compatible with the metric and orientation. The compatible complex structures on $\bR^{2m}$ are parametrised by the homogeneous space $\Sigma=SO(2m)/U(m)$ which has a standard $SO(2m)$-invariant complex structure. We get a  bundle $\pi:Z\rightarrow M$ with fibre $\Sigma$ such that a point $z$ of $\pi^{-1}(p)$ defines a compatible complex structure $I_{z}$ on $TM_{p}$. The Levi-Civita connection of $(M,g)$ defines a horizontal subbundle $H\subset TZ$, complementary to the tangent bundle along  the fibres (the vertical subbundle). We define an almost complex structure $J$ on $Z$ by specifying that at a point $z\in\pi^{-1}(p)$:
\begin{itemize}\item $J:TZ_{z}\rightarrow TZ_{z}$ preserves the horizontal and vertical subbundles.
\item $J:H_{z}\rightarrow H_{z}$ is equal to $I_{z}$ under the isomorphism  $d\pi: H_{z}\rightarrow TM_{p}$.
\item The restriction of $J$ to the vertical subspace is given by the standard complex structure on $\Sigma$. 
\end{itemize}

It is not hard to check that this almost-complex structure depends only on the conformal class of the metric $g$.

We now restrict attention to $4$-manifolds, with $m=2$. Then $\Sigma$ is just the Riemann sphere. A compatible complex structure $I$ on $\bR^{4}$ defines a 2-form $\omega$ by the usual formula $\omega(\xi,\eta)= \langle I\xi,\eta\rangle$ and this is a self-dual form of length $\sqrt{2}$. For example, if we take the standard complex structure with complex co-ordinates
$z_{1}= x_{1}+ ix_{2}, z_{2}=x_{3}+ix_{4}$ then $\omega= dx_{1}dx_{2}+dx_{3}dx_{4}$.
This gives a 1-1 correspondence between $\Sigma$ and the unit sphere in $\Lambda^{2}_{+}$.  Similarly $Z$ is identified with the unit sphere bundle in the rank $3$ real vector bundle $\Lambda^{2}_{+}\rightarrow M^{4}$. The antipodal map on the fibres defines a map $\tau:Z\rightarrow Z$ which reverses the almost-complex structure.

The curvature tensor $\Riem$ of a Riemannian manifold is a section of $\Lambda^{2}\otimes \Lambda^{2}$. In  dimensions greater than $4$ it has a decomposition into three irreducible components: the trace-free Ricci curvature $\ric_{0}$, the scalar curvature $R$ (which is the trace of the Ricci curvature) and the Weyl curvature, which is conformally invariant. 
For an oriented $4$-manifold there is another decomposition. We split $\Lambda^{2}$ into self-dual and anti self-dual parts and we get 4 components $R_{++}\in \Lambda^{2}_{+}\otimes \Lambda^{2}_{+}$ etc. of $\Riem$. Then 
\begin{itemize} \item $R_{+-}$ and $R_{-+}$ are naturally identified with the trace-free Ricci tensor ${\ric}_{0}$.
\item The traces of $R_{++}, R_{--}$ are both equal to 1/4  the scalar curvature.
\item The trace-free parts $W_{+}, W_{-}$ of $R_{++}, R_{--}$ are  the self-dual and anti-self-dual parts of the Weyl curvature.
\end{itemize}

The main foundational result then is
\begin{thm} For an oriented Riemannian $4$-manifold $(M,g)$ the  almost-complex structure on $Z$ is integrable if and only if $W_{+}=0$. 
\end{thm}

In this case the space $Z$ is called the twistor space of $(M,g)$. Atiyah, Hitchin, Singer use the opposite orientation convention, so they consider manifolds with $W_{-}=0$ called {\it self-dual 4-manifolds}, whereas we are considering anti self-dual manifolds. The $4$-manifold $M$ with its conformal structure can be recovered from $Z$ as the space of \lq\lq real lines''---embedded complex curves of genus $0$ which are preserved by the antiholomorphic involution $\tau$.  

Before discussing the proof of Theorem 2 we consider the model case when $M=\bR^{4}$. We use a different description of $\Sigma$ as $\bP(V^{+})$ where $V^{+}$ is the positive spin space. Let $\bP^{3}= \bP(V^{+}\oplus V^{-})$ so $\bP^{3}$ is a complex projective space containing two distiguished lines $\bP(V^{+})$ and $\bP(V^{-})$. Write $\bP^{3}_{*}$ for the complement of $\bP(V^{-})$. Then there is a holomorphic fibration $\varpi:\bP^{3}_{*}\rightarrow \bP(V^{+})$ with  fibres  affine complex planes whose projective completions correspond to the family of projective planes through the line $\bP(V^{-})$. 

The spin spaces $V^{+}, V^{-}$ can be regarded as $1$-dimensional quaternionic vector spaces and we have an isomorphism  $c: \bR^{4}\rightarrow {\rm Hom}_{\bH}(V^{+}, V^{-})$. For a non-zero $\psi_{+}\in V^{+}$ the map
$$x\mapsto (\psi_{+}, c(x)(\psi_{+}))$$
induces a 1-1 correspondence between $\bR^{4}$ and the affine complex plane
$\varpi^{-1}([\psi_{+}])$ and the induced complex structure on $\bR^{4}$ is  that specified by $[\psi_{+}]\in \bP(V^{+})$ and the isomorphism $\Sigma=\bP(V^{+})$. This gives an equivalence between the twistor space $Z(\bR^{4})$ and $\bP^{3}_{*}$. The twistor space of the conformal compactification $S^{4}$ is $\bC\bP^{3}$. Another way of expressing things is to start with the quaternionic vector space $\bH^{2}$. Then we can take the quaternionic projective space 
$$  S^{4}= \bH\bP^{1}= \bH^{2}\setminus \{0\}/ \bH^{*}, $$
or the complex projective space
$$  \bC\bP^{3}= \bH^{2}\setminus \{0\}/\bC^{*} $$
and we see that there is a natural map $\pi:\bC\bP^{3}\rightarrow S^{4}$ which is the twistor fibration.

We now return to the foundational Theorem 2. The Newlander-Nirenberg theorem states that an almost-complex structure is integrable if and only if its Nijenhius tensor vanishes. Given this, the proof of Theorem 2 is a differential geometric calculation which can be done in various ways. Atiyah, Hitchin and Singer work on the total space of a $\bC^{*}$-bundle over $Z$ (which is $\bH^{2}\setminus \{0\}$ in the case of $M=S^{4}$) and consider the \lq \lq twistor equation'',  an overdetermined linear equation for spinor fields on $M$. A more direct approach  exploits the conformal invariance by varying the Riemannian metric. If $p$ is a point in a Riemannian manifold we can choose a conformally equivalent metric whose Ricci curvature vanishes at $p$. This follows from the fact that the change in the Ricci curvature under a conformal factor $e^f$ is given by the Hessian of $f$ plus lower order terms. We will apply this to our $4$-manifold $M$. Another way of expressing the Newlander-Nirenberg theorem is that an almost-complex structure $J$ on a manifold $W$ is integrable if and only for each point $w\in W$ we can choose local co-ordinates centred at $w$ such that structure $J$ agrees with the standard structure up to first order--i.e. the first derivatives of the almost-complex  structure expressed in these co-ordinates vanish at the origin. We apply this to a point $z\in \pi^{-1}(p)\subset Z(M)$. Choose a metric in the conformal class such that the Ricci curvature vanishes at $p$. In the decomposition of the Riemann curvature tensor the curvature of the bundle $\Lambda^{+}$ is given by $R_{++}$ and $R_{+-}$. So in our situation, with $W_{+}+0$,  the curvature of the bundle $\Lambda^{+}$ vanishes at $p$. This means that we can choose a geodesic co-ordinate chart on $M$ and a trivialisation of $\Lambda^{+}$ around $p$ in which the connection matrix of $\Lambda^{+}$ vanishes to first order. It is then clear that the almost-complex structure on $Z(M)$ agrees with that on $Z(\bR^{4})$ to first order around $\pi^{-}(p)$ and since the latter is integrable, as we have seen above, the same holds for $Z(M)$. 


The upshot of this theory is that an (anti) self-dual conformal structure is encoded in the complex geometry of a complex 3-fold $Z$ with antiholomorphic involution $\tau$. Similarly, the ASD Yang-Mills instantons are encoded as holomorphic bundles over $Z$ by a construction originating with  Ward \cite{kn:Ward}. We consider first the noncompact structure group $GL(r,\bC)$. Recall from subsection 1.1 that a holomorphic vector bundle $E\rightarrow Z$ has an intrinsic $\db$-operator
$$  \db_{E}: \Omega^{0}(E)\rightarrow\Omega^{0,1}(E), $$
which has an extension to forms $\Omega^{p,q}(E)$ with $\db_{E}^{2}=0$ and hence the $(0,2)$ component $F^{0,2}$ of the curvature of a compatible connection $\nabla_{E}$ vanishes. A version of the Newlander-Nirenberg theorem states that conversely any connection  on a $C^{\infty}$ bundle with $F^{0,2}=0$ defines a holomorphic structure. Now let $\uE$ be a bundle with connection $\underline{\nabla}$ on the ASD $4$-manifold $M$ and let $E=\pi^{*}(\uE)$---a bundle over $Z$ with a pulled-back connection $\nabla=\pi^{*}(\underline{\nabla})$. The counterpart to Theorem 2 is
\begin{thm}
The curvature of $\nabla$ has $F^{0,2}=0$ if and only if $\underline{\nabla}$ is an ASD instanton over $M$.
\end{thm}

This is a much easier calculation. In the horizontal and vertical splitting of $TZ$ the curvature of $\nabla$ is purely horizontal, since it is lifted from the base. To understand the horizontal component we have to consider the $\pm$ self-dual decomposition for $2$-forms on $\bC^{2}$, equipped with its standard metric. The $(p,q)$ decomposition of the complexified forms is
$$  \Lambda^{2}_{\bC}= \Lambda^{2,0}\oplus \Lambda^{1,1}\oplus \Lambda^{0,2}. $$ We have a $U(2)$-invariant metric form $\omega$ which is of type $(1,1)$ so $\Lambda^{1,1}$ decomposes as the $\bC. \omega \oplus \Lambda^{1,1}_{0}$. One readily sees that the complexified anti self-dual forms are just $\Lambda^{1,1}_{0}$ and the complexified self-dual forms are $\bC.\omega\oplus\Lambda^{2,0}\oplus \Lambda^{0,2}$. This implies one direction in the statement: if $\underline{\nabla}$ is an ASD instanton the pull-back curvature has type $(1,1)$. For the converse one checks that a self-dual form on $\bR^{4}$ has non-zero $(0,2)$-component for some complex structure.

Again, one can recover the ASD connection $\underline{\nabla}$ from holomorphic
data on $Z$. By the integrability theorem, the connection $\nabla$ defines a holomorphic structure on the bundle $E$ over $Z$ and this is clearly holomorphically trivial on the fibres of $\pi$---the real lines. Conversely if $E$ is any  holomorphic bundle over $Z$ which is holomorphically trivial on the real lines then we can define a $C^{\infty}$ vector bundle $\underline{E}$ over $M$ by a \lq\lq direct image'' construction: the fibre of $\underline{E}$ at a point $p\in M$ is the space of holomorphic sections of $E$ over $\pi^{-1}(z)$. We define a connection on $\uE$ in the following fashion. Let $\underline{s}$ be a section of  $\uE\rightarrow M$, by definition this gives a section $s$ of $E$, holomorphic along the fibres. For a point $p$ in $M$ we define a subset $I_{p}$ of sections $\underline{s}$ by saying that $s\in I_{p}$ if $\db_{E} s$ vanishes along $\pi^{-1}(p)$. Then one checks that there is a unique connection $\underline{\nabla}$ on $E$ such that $\underline{s}$ is in $I_{p}$ if and only if $\underline{\nabla}\underline{s}$ vanishes at $p$. (The verification of this uses the fact that the holomorphic bundle $E$ is trivial on the first formal neighbourhood of each fibre.)

\
 The conclusion is that the study of instantons  with structure group $GL(r,\bC)$ on the ASD manifold $M$is equivalent to the study of rank $r$ holomorphic bundles on $Z$ which are trivial on all real lines. The theory can be extended to general  groups. For structure group $U(r)$ one needs the additional condition of a \lq\lq real structure'' on the bundle: an isomorphism $\overline{E}=\tau^{*}(E^{*})$.

\

\

In higher dimensions the almost-complex structure on the \lq\lq twistor space'' with fibre $SO(2m)/U(m)$ is almost never integrable. There is a good higher dimensional version of the theory, developed by Salamon  \cite{kn:SMS1},\cite{kn:SMS2} for quaternionic and  quaternion K\"ahler base manifolds. 

\subsection{The ADHM construction}

One of Atiyah's most decisive  contributions to Yang-Mills theory was his construction (in joint work with Hitchin, contemporaneous with independent work of Drinfeld and Manin) of the general solution of the Yang-Mills instanton equations over $S^{4}$ (for compact classical groups). This was preceded by some earlier work with Ward \cite{kn:MFAW} which we will touch on later. 

For the moment we ignore real structures: then the problem is to describe holomorphic bundles $E$ over $\bC\bP^{3}$ which are trivial on the generic line, so in particular $c_{1}(E)=0$. Much of the relevant theory had been developed earlier by Horrocks and Barth. The basic idea is to construct $E$ as the cohomology of a \lq\lq monad'' or 2-step complex of the form:
\begin{equation} \uU(-1)\stackrel{a}{\rightarrow} \uV\stackrel{b}{\rightarrow} \uW(1). \end{equation}
Here $U,V,W$ are complex vector spaces; $\uU,\uV,\uW$ are the corresponding trivial bundles over $\bC\bP^{3}$ and the $(p)$ notation means tensor product with the line bundle $\cO(p)$ as in 1.2. The maps $a,b$ are holomorphic bundle maps with $a$ injective, $b$ surjective and $b\circ a=0$. Such a monad defines a bundle $E$ as the cohomology ${\rm Ker}\ b/{\rm Im}\ a$. Notice that a monad description for $E$ gives one for $E^{*}$, when we replace $U,V,W$ by their duals and take the transposed maps.

\

The approach of Barth and Horrocks is based on a characterisation of sums of line bundles, due to Horrocks \cite{kn:Horrocks}.
\begin{prop}
A holomorphic vector bundle $V$ over $\bC\bP^{3}$ is a direct sum of line bundles if and only if the cohomology groups $H^{1}(V(p)), H^{2}(V(p))$ vanish for all $p$.
\end{prop}

Here one direction is completely standard. In general line bundles $\cO(p)$ over projective spaces   can have non-zero cohomology only in the top and bottom dimension. In the case of $\bC\bP^{3}$, the cohomology $H^{0}(\cO(p))$ is nonzero for $p\geq 0$ and $H^{3}(\cO(p))$ is non-zero for $p\leq 4$. Serre duality relates $p$ and $-4-p$. For $p=-1,-2,-3$ the line bundle $\cO(p)$ has no cohomology and for the analysis of the  vector bundle $E$ the key roles  are played by the cohomology groups $H^{1}(E(-1)), H^{1}(E(-2)) $ and their Serre duals $H^{2}(E^{*}(-3)), H^{2}(E^{*}(-2))$.

The main technique in this approach is  the construction of bundles as extensions, something which was prominent in Atiyah's earlier work discussed in Part 1. Given a bundle $E$ and  a vector space $W$, consider extensions
\begin{equation}  0\rightarrow E\rightarrow Q\rightarrow \uW(1)\rightarrow 0. \end{equation}
These are classified by $H^{1}(E(-1))\otimes W^{*}$. Take $W$ to be 
$H^{1}(E(-1))$, so the extensions are classified by $W\otimes W^{*}={\rm Hom}\ (W,W)$ and there is a preferred class given by the identity. This class defines a bundle $Q$ and exact sequence (21). Taking the tensor product with $\cO(-1)$ we get a long exact sequence, part of which is 
$$ W \stackrel{\partial}{\rightarrow} H^{1}(E(-1))\rightarrow H^{1}(Q(-1))\rightarrow 0 \rightarrow H^{2}(E(-1))\rightarrow H^{2}(Q(-1))\rightarrow 0 $$
The coboundary map $\partial$ is the product with the extension class, hence an isomorphism, and we see that 
$H^{1}(Q(-1))=0$ while $H^{2}(Q(-1))$ is isomorphic to $H^{2}(E(-1))$. Dually, we consider a vector space $U$ and extensions
\begin{equation}  0\rightarrow \uU(-1)\rightarrow K\rightarrow E\rightarrow
0, \end{equation}
classified by $H^{1}(E^{*}(-1))\otimes U$. We take $U$ to be the dual of $H^{1}(E^{*}(-1))$, which is $H^{2}(E(-3))$, and the extension corresponding to the identity. Now we find that $H^{1}(K(-3))=0$ while $H^{2}(K(-3))$ is isomorphic to $H^{2}(E(-3))$. In brief, passing from $E$ to $Q$ \lq\lq kills'' $H^{1}(\ \ \ (-1))$ while preserving $H^{2}(\ \ \ (-1))$ while passing from $E$ to $K$ kills $H^{1}(\ \ \ (-3))$ while preserving $H^{2}(\ \ \ (-3))$.  

Since the inclusion of $E$ in $Q$ induces an isomorphism on $H^{1}(\ \ \ (-1))$ the extension (21) prolongs to an extension $\uU(-1)\rightarrow \cV\rightarrow Q$ for some bundle $\cV$ which fits into a diagram:
\

\

\begin{picture}(1000,1000)(-300,0)
\put(300,30){$0$}
\put(500,30){$0$}
\put(100,270){$0$}
\put(300,270){$E$}
\put(500,270){$Q$}
\put(700,270){$W(1)$}
\put(990,270){$0$}
\put(100,510){$0$}
\put(300,510){$K$}
\put(500,510){${\cal V}$}
\put(700,510){$W(1)$}
\put(990,510){$0$}
\put(300,750){$U(-1)$}
\put(500,750){$U(-1)$}
\put(300,990){$0$}
\put(500,990){$0$}
\put(140,290){\vector(1,0){160}}
\put(340,290){\vector(1,0){160}}
\put(540,290){\vector(1,0){160}}
\put(820,290){\vector(1,0){160}}
\put(140,530){\vector(1,0){160}}
\put(340,530){\vector(1,0){160}}
\put(540,530){\vector(1,0){160}}
\put(820,530){\vector(1,0){160}}
\put(310,970){\vector(0,-1){160}}
\put(510,970){\vector(0,-1){160}}
\put(310,730){\vector(0,-1){160}}
\put(510,730){\vector(0,-1){160}}
\put(310,490){\vector(0,-1){160}}
\put(510,490){\vector(0,-1){160}}
\put(310,250){\vector(0,-1){160}}
\put(510,250){\vector(0,-1){160}}
\end{picture}

By construction, there are maps $a:\uU(-1)\rightarrow \cV, b:\cV\rightarrow \uW(1)$ with $ba=0$ and such that $E={\rm Ker} \ b/{\rm Im}\  a$. So the task is to show that, for a bundle $E$ corresponding to an instanton, the bundle $\cV$ is trivial. Following through long exact cohomology sequences, one finds that $H^{1}({\cal V}(-1))$ and $H^{2}({\cal V}(-3))$ vanish but
$$  H^{1}(\cV(-2))= H^{1}(E(-2))\ \ ,\ \  H^{2}(\cV(-2))= H^{2}(E(-2)). $$
So a necessary condition for $\cV$ to be trivial is that $H^{1}(E(-2))= H^{2}(E(-2))=0$.
Barth, and Barth and Hulek, \cite{kn:BH} established a converse. 
\begin{prop}
If $E$ is trivial on some line and $H^{1}(E(-2))=H^{2}(E(-2))=0$ vanish then the bundle $\cV$ is trivial, and hence $E$ arises from a monad.
\end{prop}
(Here the cohomology groups $H^{2}(E(-2))$ and $H^{1}(E^{*}(-2))$ are Serre duals, so the two vanishing conditions are interchanged by interchanging $E$ and $E^{*}$. In fact, Barth and Hulek considered bundles with $E\cong E^{*}$.) The proof is to show that, under these assumptions, all the intermediate cohomology of $\cV$ vanishes so that Horrocks' result can be applied. Then one shows that  $H^{0}(\cV(-1)), H^{0}(\cV^{*}(-1))$ vanish, which implies that  the line bundles in the sum given by Horrocks' theorem must be trivial. 

The final step then, in the proof of the ADHM result, is to show that  a bundle $E$ arising from a $U(r)$ instanton on $S^{4}$ satisfies the vanishing conditions in Proposition 6. Replacing $E$ by $E^{*}$ it suffices to show that $H^{1}(E(-2))=0$. This uses some further twistor theory, for linear equations. In general there are Penrose correspondences between cohomology classes on the twistor space and solutions of linear PDE on $S^{4}$. In fact this is a local theory, so it applies to open sets $U$ in $S^{4}$ and $\pi^{-1}(U)$ in $\bC\bP^{3}$. In the case at hand the relevant correspondence is between $H^{1}(E(-2))$ and solutions of the equation
\begin{equation}  \unabla^{*}\unabla s + \frac{R}{6} s =0 \end{equation}
for sections $s$ of $\uE$. Here $R$ is the scalar curvature of $S^{4}$. The scalar curvature enters in making the equation conformally invariant, when  considered as applying to suitable densities.  If we take $U\subset \bR^{4}$ and use the Euclidean metric the relevant equation is just $\unabla^{*}\unabla s=0$. To give an idea of this correspondence we will  consider the case when the  bundle and connection are trivial, which goes back to the origins of Penrose's theory. So we have to go from a class $\Phi$ in $H^{1}(\pi^{-1}(U), \cO(-2)) $ to a harmonic function on $U\subset \bR^{4}$. For any line $L$ in $\pi^{-1}(U)$ the cohomology group $H^{1}(L,\cO(-2))$ is isomorphic to $\bC$ and the choice of a metric in the conformal class fixes a definite isomorphism. So the restrictions of $\Psi$ to the real lines gives a function on $U$ and we want to see that this is harmonic.  We can calculate in a small neighbourhood of the origin in $\bR^{4}$ and the corresponding line $L_{0}\subset \bC\bP^{3}$. Take a standard cover of $L_{0}$ by two open sets overlapping in a neighbourhood of the equator and thicken these to open sets $\Omega_{1}, \Omega_{2}$ in $\bC\bP^{3}$. Using \v{C}ech cohomology, the class $\Phi$ can be represented by a section $\phi$ of $\cO(-2)$ over $\Omega_{1}\cap \Omega_{2}$. Recall that the twistor space of $\bR^{4}$ can be identified with the total space of the bundle
$\cO(1)\oplus \cO(1)$ over $L_{0}$. Take a co-ordinate $z$ on $L_{0}$ and standard trivialisation of $\cO(1)$ away from $z=\infty$ so we get co-ordinates
$z,\zeta_{1}, \zeta_{2}$ over $\Omega_{1}\cap \Omega_{2}$. Now it is convenient to work in the complexification, considering all lines near to $L_{0}$. These are described as sections of $\cO(1)\oplus \cO(1)$ and  depend on four complex parameters $p,q,r,s$, corresponding to the section with equations
$$  \zeta_{1}= p z+ q \ \ \ \ \ \ \  \zeta_{2}= rz + s. $$ 
Over a neighbourhood of $L_{0}$ the bundle $\cO(-2)$ can be identified with the pull back of the cotangent bundle of $L_{0}$. Thus we can write
$\phi= f(z,\zeta_{1}, \zeta_{2}) dz $
for a holomorphic function $f$ of three variables. With this set-up, the evaluation of the class on a line is defined by integrating around an equatorial contour,  so we arrive at  the function
$$  F(p,q,r,s)= \int_{\gamma} f(z, pz+q,rz +s) dz. $$
Elementary calculus shows that $F$ satisfies the equation:
$$  \frac{\partial^{2} F}{\partial p\partial s} - \frac{\partial^{2} F}{\partial q\partial r}=0. $$
This is the complexified Laplace equation corresponding to the quadratic form $ps-qr$ on $\bC^{4}$. Making a linear change of co-ordinates to put this quadratic form into the standard shape $\sum x_{i}^{2}$ and then restricting to the real  points we get an integral representation of harmonic functions, which in fact goes back to Bateman \cite{kn:Bateman} in 1904. 

What we have discussed above is just the easiest part of the Penrose correspondence. But, assuming that correspondence,  the vanishing of $H^{1}(E(-2))$ for a bundle $E$ corresponding to a unitary instanton is immediate because 
$\unabla^{*}\unabla  + \frac{R}{6} $ is a positive operator, so the only global solution of (23) over $S^{4}$ is $s=0$.

\subsection{Topology and geometry of moduli spaces}
\subsubsection{Riemann surfaces}
In this subsection we discuss the long and influential paper {\it The Yang-Mills equations over Riemann surfaces} \cite{kn:AB1} by Atiyah and Bott.
This takes up some of the themes in Atiyah's earlier work, discussed in Section 1, and subsequent developments outlined in 1.4. Let $X$ be a compact Riemann surface of genus $g\geq 1$. For each $r$ and $0\leq d<r$ there is a moduli space $\cM_{r,d}$ of stable holomorphic vector bundles over $X$ of rank $r$ and degree $d$ and fixed determinant, and these have a Narasimhan-Seshadri description in terms of projective unitary representations. We will concentrate on bundles of rank $r=2$. The most straightforward case in some ways is to take $d=0$ when one can fix the determinant to be the trivial line bundle, but there is a major complication that there are then also reducible representations. The moduli space of irreducible representations is not compact and the obvious compactification, adjoining the reducible representations, is singular. So we will consider the  case with $r=2, d=1$ and we get a compact moduli space $\cM= \cM_{2,1}$  which is a complex manifold of dimension $3g-3$. The problem is to describe the cohomology of $\cM$. This question had been treated before from several points of view, beginning with Newstead who found the Betti numbers, and we will return to that below.

The novelty in the work of Atiyah and Bott was to bring in ideas from Yang-Mills theory. Let $Y$ be any compact Riemannian manifold and $E\rightarrow Y$ a $U(r)$ bundle. Then we have an infinite-dimensional space $\cA$ of all connections on $E$ and an infinite-dimensional group $\cG$ of bundle automorphisms which acts on $\cA$. The Yang-Mills functional $\cL(A)= \Vert F(A)\Vert^{2}$ is a $\cG$-invariant functional on $\cA$ and its critical points are the Yang-Mills connections on $E$. If we have a compact  hence finite dimensional, manifold $B$ and a function $f$ on $B$ there are well-established theories  which relate the critical points of $f$ to the topology of $B$. There are also well-established extensions of these ideas to infinite dimensions, notably to the energy functional on paths in a manifold whose critical points are geodesics. To apply these ideas to the Yang-Mills case two basic points arise.
\begin{itemize}
\item As in any variational problem, there are  fundamental analytical and PDE questions which have to be tackled and these depend very much on the dimension of the base manifold $Y$. We will say more about  this later. 
\item The action of the gauge group $\cG$ means that the critical points are never isolated and one would like to work on the quotient space $\cB=\cA/\cG$ but then there are difficulties caused by the existence of reducible connections.
\end{itemize}

To set up the picture in more detail, fix a compatible Riemannian metric on the Riemann surface $X$ so that $X$ has area $2\pi$ and fix a $C^{\infty}$ rank $2$ bundle $E\rightarrow X$ with $c_{1}(E)=1$. Fix a connection on the determinant bundle $\Lambda=\Lambda^{2}E$
with curvature $-i$ times the area form on $X$. Let $\cA$ be the space of connections on $E$ which induce the given connection on $\Lambda^{2}E$ and let $\cG$ be the group of automorphisms of $E$ which act with determinant $1$ on the fibres. The Yang-Mills connections are easy to describe. Recall that the equations are $d_{A}^{*}F(A)=0$, or equivalently $d_{A}(* F(A))=0$. In this case $*F(A)$ is a section of the bundle $\rad E\subset {\rm End} E$ and the equations say that it is covariant constant. For convenience, set $\bF(A)= - i *F(A)$, so $\bF(A)$ is a Hermitian endomorphism of $E$ with trace $1$. If $\bF(A)$ is covariant constant---{\it i.e.} $d_{A}\bF(A)=0$--- there are two possibilities:
\begin{enumerate}\item $\bF(A)= \frac{1}{2} 1_{E}$;
\item $\bF(A)$ has constant eigenvalues $(\mu,1-\mu)$ for some integer $\mu\geq 1$.
\end{enumerate}
In the first case the induced connection on the trace free part of $\rad E$ is a flat $SO(3)=PU(2)$ connection and it is straightforward to show that these precisely correspond to the projective unitary representations given by the Narasimhan-Seshadri theorem. It is also easy to show that these are absolute minimisers of the Yang-Mills functional on $\cA$. So in $\cB=\cA/\cG$ the minimum of the functional can be identified with the moduli space $\cM$ we want to study. 

 In the second case,  the eigenspaces give a decomposition $E=L\oplus L^{*}\otimes \Lambda$ where $L$ has degree $\mu$ and the connection is given by a sum of connections with constant curvature on the two line bundles. Since the connection on $\Lambda$ is fixed this is determined by a connection with constant curvature on $L$. Conversely for any connection on $L$ we can build a Yang-Mills connection on $E$. Two connections with constant curvature on $L$ differ by a flat connection, or in other words a representation $\pi_{1}(\Sigma)\rightarrow S^{1}$. These representations are parametrised by the $2g$-dimensional  torus
 $J= H^{1}(X;\bR)/H^{1}(X,\bZ)$.
The Yang-Mills functional takes the value $q(\mu)=\mu^{2}+(1-\mu)^{2}$ on these reducible connections. The conclusion is that the non-minimal critical points of the functional in $\cA/\cG$ consist of copies $J_{\mu}$ of $J$ for $\mu\geq 1$, with critical value $q(\mu)$.

The discussion so far has been based on the Riemannian geometry of $X$. Now we pass to the complex geometry. Any connection $A\in \cA$ defines a holomorphic structure on $E$ through its $\db$-operator $\db_{A}$. For line bundles this gives the classical identification of the real torus $J$ with the Jacobian of $X$,  parametrising holomorphic line bundles of degree $0$. The situation can be described by introducing another group $\cG^{c}$, consisting of automorphisms of $E$ of determinant $1$ but not necessarily unitary. Then $\cG^{c}$ acts on $\cA$ in the following way. Thinking of a connection as a covariant derivative
$$  \nabla_{A}= \db_{A}+\partial_{A} $$
an automorphism $g\in \cG^{c}$ acts by
\begin{equation}  \nabla_{g(A)}= g\circ \db_{A}g^{-1} + (g*)^{-1}\partial_{A}  g^{*} . \end{equation}
Said in another way, a connection on $E$ is determined by its $\db$-operator and the given Hermitian structure. (This is  what is often called the Chern connection in complex differential geometry, see Chapter 0 in \cite{kn:GH} for example.) So we can identify $\cA$ with the space of $\db$-operators and there is an obvious action of $\cG$ on that space (which corresponds to the first term on the right hand side of (24)).

 The conclusion from the above is that the { set} of equivalence classes of rank-2 holomorphic vector bundles over $X$ with fixed determinant $\Lambda$ can be identifies with the set of orbits of $\cG^{c}$ in $\cA$. We write  {\it set} rather than {\it space} because the natural topology on the quotient $\cA/\cG^{c}$ is not Hausdorff, due to jumping phenomena of the kind we discussed in 1.3. This is the point of the notion of a \lq\lq stable'' bundle. Recall that in 1.2 we introduced  an invariant  $\mu(V)$ of a rank $2$ holomorphic vector bundle over $X$: the maximal degree of a line subbundle. In our situation, when $V$ has degree $1$ the stability condition is just that $\mu(V)< 1$. Then we get a $\cG^{c}$-invariant stratification of $\cA$ consisting of
\begin{itemize} \item An open set $\cA_{s}\subset \cA$ of \lq\lq stable points'', i.e. connections which define stable holomorphic structures.
\item For each $\mu\geq 1$ a stratum $C_{\mu}\subset \cA$ of connections which define a holomorphic structure $V$ with $\mu(V)=\mu$.
\end{itemize}

Let us now go back to the variational point of view. Following standard strategy and the analogy with finite dimensional problems we could consider the gradient flow of the Yang-Mills functional, which is the evolution equation in $\cA$
\begin{equation}   \frac{\partial A}{\partial t} = -d^{*}_{A} F(A), \end{equation} 
and hope to show that this converges as $t\rightarrow \infty$ to a critical point. This would give a $\cG$-invariant stratification of $\cA$, with strata labelled by the critical sets, which would be the basis for a Morse-theoretic analysis of the topology. But all this  would depend on detailed analysis of the equation (25). Atiyah and Bott realised that the complex geometry allows one to bypass this (then conjectural) analysis because the putative \lq\lq Morse'' stratification is exactly the same as the stratification above coming from the complex geometry. In more detail:
\begin{itemize} \item By definition of the moduli space of stable holomorphic bundles and the Narasimhan-Seshadri theorem there is a map $R:\cA_{s}\rightarrow \cM$;
\item Let $A$ be a point in $C_{\mu}$, so it defines a holomorphic structure $V$ which is an extension
\begin{equation} 0 \rightarrow L\rightarrow V\rightarrow \Lambda\otimes L^{*}\rightarrow 0\end{equation}
for a holomorphic line bundle $L$ of degree $\mu\geq 1$. It is easy to see that this subbundle is unique. Then we get a map $R_{\mu}: C_{\mu}\rightarrow J_{\mu}$ taking $A$ to the corresponding reducible connection on $L\oplus \Lambda\otimes L^{*}$. 
  \end{itemize}

 Then one can show from simple calculations with the flow equation (25) that {\it if} the flow exists for all time and converges to a critical point the limit points are given by the maps $R$ and $ R_{\mu}$. The essential point is that the flow preserves the $\cG^{c}$-orbits. The whole picture extends to bundles of higher rank. An unstable bundle $V$ has a canonical \lq\lq Harder-Narasimhan filtration'' so that successive quotients are semi-stable and this can be refined to a filtration
$$  0=V_{0}\subset V_{1}\dots\subset V_{k}\subset V$$ with
$ V_{i+1}/V_{i}$ stable. The gradient flow starting with a holomorphic structure of type $V$ has a limit as $t\rightarrow \infty$ corresponding to the holomorphic structure $\bigoplus V_{i+1}/V_{i}$, with a projectively flat connection on each factor. 

We now move to the second point: taking account of the $\cG$-action on $\cA$. For this, Atiyah and Bott introduced the use of equivariant cohomology in Morse Theory. For any group $\Gamma$ and space $Z$ with $\Gamma$ action the $\Gamma$-equivariant cohomology
$H^{*}_{\Gamma}(Z)$ is defined to be the ordinary cohomology of the space
$Z_{\Gamma}= Z\times_{\Gamma} E\Gamma$ where $E\Gamma$ is the total space of the universal bundle $E\Gamma\rightarrow B\Gamma$. Thus $Z_{\Gamma}$ is a fibre bundle over $B\Gamma$ with fibre $Z$. Restriction to the fibre gives a map from $H^{*}_{\Gamma}(Z)$ to $H^{*}(Z)$ and pulling back from the base gives a map from $H^{*}(B\Gamma)$ to $H^{*}_{\Gamma}(Z)$ which makes $H^{*}_{\Gamma}(Z)$ a module over the ring $H^{*}(B\Gamma)$. In the situation at hand we want to consider the equivariant cohomology $H^{*}_{\cG} (\cA)$ and relate this to the stratification. 

The general picture in equivariant theory is, very roughly, that a point in $Z$ with stabiliser
$\Gamma'\subset \Gamma$ should contribute $H^{*}(B\Gamma')$ to the cohomology.
In our case the stabiliser  in $\cG$ of a point in $J_{\mu}$ is a copy of $S^{1}$ and $BS^{1}=\bC\bP^{\infty}$ with cohomology the polynomial ring in one variable of degree $2$. Slightly more precisely, the whole critical set $J_{\mu}$ can potentially contribute $H^{*}(J_{\mu}\times \bC\bP^{\infty})$ to the equivariant cohomology of $\cA$ but this is shifted in degree by 
the index: the dimension of the negative subspace for the Hessian of the functional $\cL$. That is, from standard arguments using the existence of the stratification and the retraction maps $R, R_{\mu}$, the $\cG$-equivariant cohomology of $\cA$ is no larger than the direct sum of the ordinary cohomology
$H^{*}(\cM)$ (which we are trying to compute) with the sum over $\mu$ of the shifted $H^{*}(J_{\mu}\times \bC\bP^{\infty})$. A fundamental observation of Atiyah and Bott is the \lq\lq self-completing principle'' which states that the equivariant cohomology of the whole space $\cA$ is equal to this direct sum. Said in another way the \lq\lq Morse inequalities'' in this situation become equalities.  This follows from general properties of the equivariant set-up. The index of a point of $J_{\mu}$ can be computed to be $2g+4\mu$ and the conclusion can be expressed in terms of the equivariant Poincar\'e series
$ P(t)= \sum {\rm dim}\ H^{q}_{\cG}(\cA)$ as 
 
\begin{equation} P(t)= P_{\cM}(t)+
 \sum_{\mu=1}^{\infty} t^{2g+4\mu} \frac{(1+t)^{2g}}{(1-t^{2})}\end{equation}
where $P_{\cM}$ is the ordinary Poincar\'e polynomial of the moduli space $\cM$.

The next question is to identify the  equivariant cohomology of  the space of connections $\cA$. Since $\cA$ is contractible this is the same as the cohomology of the classifying space $B\cG$ which can be found by standard algebraic topology arguments. The answer is as follows. There is a universal bundle, which is a $\cG$-equivariant $SO(3)=PU(2)$ bundle over $X\times \cA$. Here $\cG$ acts trivially on $X$. This has a Pontrayagin class in the $4$-dimensional equivariant cohomlogy of $X\times \cA$ whicb is $H^{*}(X)\otimes H^{*}_{\cG}(\cA)$.
 Taking the K\"unneth components of this we get  classes $u\in H_{\cG}^{4}(\cA),
 \omega\in H_{\cG}^{2}(\cA)$ and a map $\nu:H_{1}(\Sigma)\rightarrow H_{\cG}^{3}(\cA)$. Then the result is that $H_{\cG}{*}(\cA)$ is freely generated by these classes, {\it i.e.}
$$  H_{\cG}^{*}(\cA)= \bQ[u,\omega]\otimes \Lambda^{*} H_{1}(X)$$
where $u$ has degree $4$, $\omega$ has degree $2$ and $H_{1}(X)$ has degree $3$. Thus the Poincar\'e series is
$$  P(t)= \frac{(1+t^{3})^{g}}{(1-t^{4})(1-t^{2})}. $$ 
and combining this with (27) gives a formula for the Poincar\'e polynomial of $\cM$:
$$  P_{\cM}(t)= \frac{1}{(1-t^{2})(1-t^{4})}\left( (1+t^{3})^{2g}- t^{2g} (1+t)^{g}\right). $$

This can be expressed in various other ways. Examining the Atiyah-Bott argument gives a refinement to describe the homology groups of $\cM$ as  representations  of the mapping class group of $X$ (which acts on $\cM$ via its description through projective flat connections). These representations are all built out of the  exterior powers $\Lambda^{j}= \Lambda^{j} H^{1}(X;\bQ)$. It is neatest to normalise the gradings about the middle dimensions, so we write
$$  \uLambda^{j}= \Lambda^{g+j} \ \ \  \uH_{i}= H_{3g-3+ i}(M;\bQ). $$
The one finds \cite{kn:SKD3} the genus-independent formula
$$  \uH_{i}= \sum a_{ij} \uLambda_{j}, $$

where $a_{ij}$ is the co-efficient of $t^{i}$ in 
$$ \frac{(t^{2j}-t^{-2j})(t^{j}-t^{-j})}{(t^{2}-t^{-2})(t-t^{-1})}. $$
For example when $g=2$ the rational homology of $\cM$---the intersection of two quadrics as described in 1.4---is $\bQ$ in dimensions $0,2,4,6$ and a copy of $H^{1}(X;\bQ)$ in dimension $3$.

\subsubsection{Other results and developments}

The arguments of Atiyah and Bott applied to bundles of all ranks and they obtained inductive formulae for the Betti numbers of moduli spaces $\cM(r,d)$ of stable bundles for $r,d$ coprime. In the simplest case considered above of $\cM(2,1)$ the Betti numbers had been found first by Newstead \cite{kn:Newstead1}. General formulae had been found by Harder and Narasimhan \cite{kn:HN} using the Weil conjectures and counting points in moduli spaces defined over finite fields. Atiyah and Bott pointed out intriguing parallels between the different approaches---one related to physics and one to number theory--- which were developed further by Asok, Doran and Kirwan \cite{kn:ADK} and, more recently, Gaitsgory and Lurie \cite{kn:GL}.

Atiyah and Bott also found further information about the topology of these moduli spaces. They showed that they are simply connected and the integral homology has no torsion and they obtained  results about the cohomology ring. Returning again, for simplicity, to the case of $r=2,d=1$ their arguments showed that  the restriction map from $H^{*}(B\cG)$ to $H^{*}(\cM)$ is surjective;  we have explicit generators $\omega, u, \nu( H_{1}(\Sigma))$ for the former so it is question of finding the relations between these in $H^{*}(\cM)$. For example when $g=2$ one easily finds that $\omega^{2}= 2 u$ in $H^{4}(M)$. They also made progress on a conjecture of Newstead regarding the characteristic classes of $M$, which they showed reduced to the conjecture $p_{1}^{g}=0$ where $p_{1}$ is the first Pontrayagin class of the tangent bundle

The decade following the Atiyah and Bott paper saw many further developments in the understanding of the topology of these moduli spaces of which we just mention a few. 
 \begin{itemize}\item King and Newstead gave a complete description of the relations in the cohomology rings $H^{*}(\cM)$  for rank 2 bundles \cite{kn:KN}. The  higher rank case was treated by Earl and Kirwan \cite{kn:EK}.
\item Kirwan proved the Newstead conjecture $p_{1}^{g}=0$ \cite{kn:FCK2}.
\item In another direction, one can consider the integrals
\begin{equation}  \int_{M} \omega^{p} u^{q}, \end{equation}
where $p+ 2q=3g-3$ and the dimensions of the vector spaces
$$  H^{0}(M,L^{k})$$ where $L$ is the  positive line bundle over the moduli space with $c_{1}(L)=\omega$. These are related by the Riemann-Roch formula.  The vector spaces came into prominence as \lq\lq conformal blocks'' in conformal field theory, fitting later into Witten's quantum field theory interpretation, and generalisation, of the Jones polynomial. This was a major interest of Atiyah's in the late 1980's and one example of the  general notion of a Topological Qunatum Field Theory, as discussed in Dan Freed's article in this volume. In \cite{kn:MT}, Thaddeus used the Verlinde formulae for the dimensions of these vector spaces to evaluate the integrals (28).  About the same time, Witten \cite{kn:Witten}, developed a related theory, focused on  the volume of the moduli spaces, which could also be used to find the integrals. Later, Jeffrey and Kirwan \cite{kn:JK} evaluated the corresponding integrals over moduli spaces of bundles of all ranks,   using a localisation formula related to what we will discuss in 2.4.5 below. 

\item Another huge development involved moduli spaces of \lq\lq pairs'' consisting of a holomorphic bundle $E$ with extra structure. Prominent examples, which we will touch on later, are
vortices $(E,\psi)$, studied by Bradlow, Garcia-Prada and many others,  where $\psi$ is a holomorphic section of $E$, and the Higgs bundles $(E,\phi)$ introduced by Hitchin \cite{kn:NJH2}, where $\phi$ is a holomorphic $1$-form with values in ${\rm End} E$. 
\end{itemize}
\subsubsection{Moment maps: symplectic, K\"ahler and hyperk\"ahler}

One of the most influential results in the paper of Atiyah and Bott---a digression from the main theme of the paper---was their interpretation of curvature as a moment map. This is a concept from symplectic geometry and Atiyah wrote several papers about aspects of symplectic geometry in the mid-1980's. 

Recall the Hamiltonian construction in symplectic geometry. If $(V,\omega)$ is a symplectic manifold  and $H$ is a smooth function on $V$ one defines a vector field $X_{H}$ by the condition that its contraction with $\omega$ is the $1$-form $dH$. This has the property that the Lie derivative of $\omega$ along $X_{H}$ vanishes, so the flow generated by $X_{H}$ is a $1$-parameter family of symplectomorphisms,  preserving $\omega$. Conversely, at least locally, any such flow is generated by a Hamiltonian function $H$. Now suppose that a Lie group $G$ acts on $V$ preserving $\omega$. A map
$$  \mu:V\rightarrow {\rm Lie}(G)^{*}$$

is called a moment map for the action if for each $\xi$ in ${\rm Lie}(G)$ the pairing $\langle \mu,\xi\rangle$ is a Hamiltonian function for the $1$-parameter subgroup generated by $\xi$. In other words, at each point $x\in V$ the derivative $d\mu: TV_{v}\rightarrow {\rm Lie} (G)^{*}$ is the transpose of the infinitesimal action ${\rm Lie}(G)\rightarrow TV_{v}$ under the identification of $TV_{v}$ with its dual furnished by $\omega$. The   moment map is equivariant if it intertwines the actions of $G$ on $V$ and ${\rm Lie}(G)^{*}$. 

Now go back to the infinite dimensional space $\cA$ of connections on a bundle over a compact Riemann surface $X$. This is an affine space with associated vector space $\Omega^{1}(X; \rad E)$. There is a symplectic form
\begin{equation}  \Omega(a,b)= \int_{X} {\rm Tr}(a\wedge b) . \end{equation}

Note that this does not use the complex structure on $X$, only the orientation. This is a generalisation of the skew pairing on ordinary $1$-forms defined by wedge product and integration. 

The gauge group $\cG$ acts on $\cA$ preserving the symplectic form. The Lie algebra of $\cG$ is the space $\Omega^{0}(X,\rad E)$ and integration gives a map from $\Omega^{2}(X,\rad E)$ to the dual of the Lie algebra. The observation of Atiyah and Bott is that the map from $\cA$ to $\Omega^{2}(X,\rad E)$ given  by the curvature of a connection is an equivariant moment map for the action of $\cG$. This paved the way to a number of important subsequent developments. In one direction, in symplectic geometry, with an equivariant moment map $\mu$ as above,  one considers the \lq\lq symplectic'' or \lq\lq Marsden-Weinstein'' quotient
$$  V//G= \mu^{-1}(0)/G$$
  which has (ignoring possible singularities) an induced symplectic form. Applied to the space of connections, this gives a symplectic form on the moduli space of irreducible  flat connections. For a bundle $E$ with nonzero degree there is a slight modification of this. In the general picture, with a group $G$ acting on $(V,\omega)$,   if the Lie algebra of $G$ has a non-trivial centre a moment map is not unique since one can add any element of the dual of the centre. The basic example is when $G=S^{1}$ so a  
 moment map is a Hamiltonian function and one can add a constant. For a bundle with non-zero degree we modify the moment map so that the zeros are given by connections with constant central curvature and we get a symplectic form on the moduli spaces we considered before. About the same time, Goldman defined symplectic structures on moduli spaces of representations of the fundamental group of a surface into any semi-simple Lie group \cite{kn:WMG}. Taking the group $SL(2,\bR)$, he made a connection with results of Wolpert \cite{kn:SW} on the symplectic geometry of Teichmuller space.

In another direction, the identification of the curvature as a moment map fits the Narasimhan-Seshadri theorem into a general theme of \lq\lq equality of symplectic and complex quotients'' which emerged at around the same time as the Atiyah-Bott paper. This is a theory which initially applies in finite dimensional situations. Suppose that $V$ is a complex manifold, $\omega$ is a K\"ahler form and $G$ is a compact Lie group acting on $V$ and preserving  both the complex structure and $\omega$. There is a complexification $G^{c}$ and the action extends to a holomorphic action of $G^{c}$ on $V$. The full quotient $V/G^{c}$ will  not usually be a well-behaved object--it will not be Hausdorff--- but under suitable technical hypotheses one defines a set of \lq\lq stable'' points  $V^{s}\subset V$ such that the quotient $V^{s}/G^{c}$ is Hausdorff. Then the general principle is that this complex quotient agrees with the symplectic quotient defined via  a moment map. Further, the pictures are compatible in that the induced symplectic form on the quotient space is a K\"ahler form with respect to the visible complex structure in the complex description. In the case of moduli spaces of flat connections this gives a simple conceptual setting for the K\"ahler structures which had been written down before by Narasimhan \cite{kn:MSN}.

This principle of equality of symplectic and complex quotients sheds light on the discussion in 1.3 above of small resolutions.
Take $G=S^{1}$ acting on $\bC^{2}\times \bC^{2}$ with weight $+1$ on the first factor and $-1$ on the second. With the standard symplectic form a Hamiltonian is
$$  H(x,y)= \vert x\vert^{2}-\vert y\vert^{2}. $$
We are free to add a constant so for each $\epsilon\in \bR$ we have a symplectic quotient $H^{-1}(\epsilon)/S^{1}$. Taking $\epsilon>0$ means that we are considering pairs $(x,y)$ with $x\neq 0$ and the symplectic quotient agrees with the quotient $U^{+}$ that we considered before. Similarly, taking $\epsilon=0$ or $\epsilon<0$ gives the other two quotient spaces that we considered. In other words, the \lq\lq Atiyah flop'' interchanging the two small resolutions appears  as  the variation of the symplectic quotient with respect to the level set of the Hamiltonian. In the complex geometry picture this enters through the fact that one has to make a choice of stability condition. The study of variation of moduli spaces under change of stability conditions has become a huge area. As one example, we  mention the work of Thaddeus \cite{kn:MT2} in which he considered vortices---pairs $(E,\psi)$ consisting of a rank 2 bundle $E$ and a holomorphic section $\psi$. Thaddeus takes bundles with odd degree and  fixed determinant.  There is a family of stability conditions depending on a real parameter $\sigma$ and hence a family of moduli spaces $N_{\sigma}$. When $\sigma$ is very large the condition just requires that the bundle $E$ is stable, so in this regime $N_{\sigma}$ fibres over the moduli space $\cM_{2,1}$ with fibre $\bP(H^{0}(E))$ over $E$. When $\sigma$ is small the stability condition just requires that $\psi$ has no zeros and $N_{\sigma}$ is a projective space parametrising extensions up to scale. Between these extremes there are a finite number of values of $\sigma$ where the moduli space changes, as in the model example above. Thaddeus was able to describe these changes explicitly and, using this, was able to give a proof of the Verlinde formula and also another approach to the homology of the moduli spaces.

The equivariant Morse theory techniques of Atiyah and Bott were developed in the thesis of Kirwan, supervised by Atiyah, obtaining a general machine for calculating cohomology of quotient spaces in symplectic and algebraic geometry \cite{kn:FCK1}. As one of the  simplest  examples, Kirwan's theory calculates the Betti numbers of the quotient of the set of configurations  of $d$ points in the Riemann sphere by the action of the the M\"obius transformations (assuming, for technical reasons, that $d$ is odd).

In infinite-dimensional situations, such as the action of $\cG$ on $\cA$, this principle of equality of symplectic and complex quotients is a guide rather than  a general theorem. But the idea has been important in subsequent developments in complex differential geometry. These include examples related to gauge theory, some of which we will discuss further below, and also to other structures such as K\"ahler metrics. The author has written about this in many places so  we will just refer to \cite{kn:SKD4}, \cite{kn:SKD5}.

These moment map ideas are also important in four-dimensional gauge theory and mesh in with  other work of Atiyah.  In abstract, suppose now that $V$ is a {\it hyperk\"ahler manifold}. That is, $V$ has a Riemannian metric $g$, three compatible complex structures $I_{1}, I_{2}, I_{3}$ giving an action of the quaternions and three corresponding K\"ahler forms $\omega_{1}, \omega_{2}, \omega_{3}$. Suppose that $G$ acts on $V$ preserving all this structure. Then a hyperk\"ahler moment map is a map 
$$  \underline{\mu}: V\rightarrow {\rm Lie}(G)^{*}\otimes \bR^{3}$$
whose components give moment maps as considered before for the three symplectic forms $\omega_{i}$. This set-up was considered by Hitchin {\it et al.} \cite{kn:HKLR}, who showed that, under suitable technical hypotheses, the quotient
$$ V///G= \underline{\mu}^{-1}(0)/G, $$ 
has an induced hyperk\"ahler structure. To explain this in outline,  write $\underline{\mu}=(\mu_{1},\mu_{2}, \mu_{3})$ and fix one complex structure, say $I_{1}$. Then one finds that $\mu_{2}+i\mu_{3}$ is holomorphic with respect to $I_{1}$, so its zero set is an $I_{1}$-complex subvariety. Then (ignoring possible singularities) the quotient $V///G$ can be regarded as the K\"ahler quotient of this complex
subvariety, as considered before. Replacing $I_{1}$ by $I_{2}$ and $I_{3}$ one sees the three complex structures on $V///G$. 

This theory can be applied in two ways to the moduli spaces of instantons over $\bR^{4}$, which end up giving the same hyperk\"ahler structures. For the first. we work in an infinite dimensional setting, as in the paper of Atiyah and Bott.  Let $\omega_{1}, \omega_{2}, \omega_{3}$ be the standard basis for the self-dual $2$-forms on $\bR^{4}$. Then for tangent vectors $a,b$ to the space $\cA$ of connections on $E\rightarrow \bR^{4}$ we have symplectic forms
$$\Omega_{i}(a,b)= \int_{\bR^{4}} {\rm Tr}(a\wedge b) \wedge \omega_{i}, $$
making $\cA$ formally into an infinite dimensional hyperk\"ahler manifold
and the self-dual part of the curvature $F^{+}(A)$ has three components $F^{+}_{i}(A)$.
The basic fact is that the $F^{+}_{i}$ are moment maps for the gauge group action on $\cA$ with respect to the forms $\Omega_{i}$: in other words, $F^{+}$ is a hyperk\"ahler moment map. Since $\bR^{4}$ is not compact one should really be more careful in specifying what kind of behaviour at infinity is allowed but suffice it here to say that the relevant gauge group should be the \lq\lq based gauge group'' of automorphisms which fix the fibre over infinity. Thus the quotient is the moduli space of instantons \lq\lq framed'' at infinity.  The conclusion then is that these  moduli spaces of framed instantons have hyperk\"ahler structures given by the quotient construction. 

 The other approach is more elementary and uses the ADHM description. Recall that this identifies moduli spaces of instantons with equivalence classes of monads, which are given by matrix data. We write down explicitly what this data is, after suitable normalisations. To do this it is convenient to single out one complex structure on $\bR^{4}$. Then for an instanton with structure group $U(r)$ and Chern class $k$ one finds that
 the matrix data consists of 
linear maps $$\alpha_{1},\alpha_{2}:\bC^{k}\rightarrow \bC^{k},
 P:\bC^{r}\rightarrow \bC^{k}, Q:\bC^{k}\rightarrow \bC^{r}$$
satisfying the equations
\begin{equation}  [\alpha_{1}, \alpha_{2}]+ PQ=0 , \end{equation}
\begin{equation} [\alpha_{1},\alpha_{1}^{*}]+ [\alpha_{2}, \alpha_{2}^{*}]+ PP^{*}- Q^{*}Q=0 , \end{equation}
and also certain nondegeneracy conditions. The framed moduli space is given by the quotient of such solutions $(\alpha_{1}, \alpha_{2}, P,Q)$ by the obvious action of $U(k)$. (The obvious action of $U(r)$ changes the framing.) Now the point is that the set of all quadruples of  maps $(\alpha_{1}, \alpha_{2}, P,Q)$ is a quaternionic vector space:
  $$ {\rm Hom}(\bC^{k},\bC^{k})\otimes_{\bC} \bH \oplus {\rm Hom}(\bC^{k}, \bC^{r})\otimes_{\bC}\bH, $$
and the quadratic expressions on the left hand sides of (30),(31) make up the hyperk\"ahler moment map for the action of $U(k)$. Thus the moduli space appears as a finite dimensional hyperk\"ahler quotient. Geometrically, the choice of complex structure gives a projective plane $\bC\bP^{2}$ inside the twistor space and the the equation (30) is the condition that the data $(\alpha_{1},\alpha_{2},P,Q)$ provides a monad defining a holomorphic vector bundle over this plane, which is just the restriction of the bundle over $\bC\bP^{3}$ considered in the twistor approach. In fact this construction sets up a 1-1 correspondence between the framed moduli spaces of instantons and moduli spaces of holomorphic bundles over $\bC\bP^{2}$ trivialised over the line at infinity. This was shown in \cite{kn:SKD1}, confirming a suggestion of Atiyah in Chapter VII of \cite{kn:MFA5}.

At this point we  make a digression which will connect again with Atiyah's early work on double points. As we said in 1.3,  the remarkable feature is that a 2-dimensional double point singularity can be resolved in a family (after taking a covering). This leads to a famous story in complex geometry, through work of Brieskorn and others, involving \lq\lq ADE singularities''. The finite subgroups of $SU(2)$ are in 1-1 correspondence with the Lie algebras of type A,D,E. For each such group $\Gamma$ the quotient
$\bC^{2}/\Gamma$ has the special property that its resolution and smoothings are diffeomorphic and it can be resolved in a family. The intersection pattern of the exceptional curves in the resolution yields the Dynkin diagram of the Lie group. The ordinary double point singularity is the case $A_{1}$ where the subgroup $\Gamma$ is $\pm 1$ and there is just one exceptional curve. Moving now to differential geometry, Gibbons and Hawking introduced the notion of \lq\lq ALE gravitational instantons''. These are hyperk\"ahler 4-manifolds which are asymptotically locally Euclidean {\it i.e.}  the structure asymptotic at infinity to the flat structure on $\bC^{2}/\Gamma$ for some $\Gamma\subset SU(2)$ as above. In the simplest case when $\Gamma=\{\pm 1\}$ such a metric had
been found earlier, independently,  by Calabi and Eguchi-Hanson. For the other $A_{k}$ singularities, Gibbons and Hawking gave an explicit construction based on harmonic functions on $\bR^{3}$ with poles. Hyperk\"ahler structures on $4$-manifolds are (anti) self-dual so have   twistor spaces. Hitchin found the twistor description of the Gibbons-Hawking metrics and showed that it fitted beautifully with the Brieskorn theory.  Such a metric has a family of complex structures parametrised by the $2$-sphere. For two special points on the sphere the complex structure is that of the resolution of the singularity and for the other points of the smoothing. In his paper \cite{kn:MFA7} Atiyah used Hitchin's twistor description to find explicit formulae for the Green's function on these $A_{k}$ gravitational instantons. Regarding a pair of points in the 4-manifold as a pair of lines in twistor space, he showed that the Green's function corresponds to the Serre class in an ${\rm Ext}$ group, related also to a paper of Atiyah and Ward \cite{kn:AW} on the construction of bundles on twistor spaces. 

To get back to hyperk\"ahler quotients; the question left open after the developments outlined above was the existence of gravitational instantons corresponding to the D and E families of Lie groups. This was answered decisively by Kronheimer in his thesis supervised by Atiyah \cite{kn:Kronheimer}. Kronheimer used a quotient construction which has some similarities to the ADHM description of instanton moduli spaces. For any finite subgroup $\Gamma\subset SU(2)$, Kronheimer considers the regular representation $R$. Then $R\otimes \bC^{2}$ is a quaternionic representation of $\Gamma$. Let $M= \left(R\otimes \bC^{2}\right)^{\Gamma}$ be the space of vectors fixed by $\Gamma$ and let $G$ be the subgroup of $U(R)$ of transformations which commute with $\Gamma$.  Then $G$ acts on $M$ preserving a hyperk\"ahler structure and Kronheimer showed that the hyperk\"ahler quotient is the desired ALE space. One has a family of these quotients obtained by varying the level set of the moment map and in this way Kronheimer  obtained a complete classification of all ALE gravitational instantons.

\

These ideas also apply to solution of the Bogomolny equation or {\it monopoles} which were a major research interest of Atiyah. The Bogomolony equation is an equation on $\bR^{3}$ obtained from the dimension reduction of the instanton equation on $\bR^{4}$. One component of the connection becomes a Higgs field: so the data is a connection $A$ on a bundle $E$ over $\bR^{3}$ and a section $\phi$ of $\rad E$. The Bogomolny monopole equation is
\begin{equation}  d_{A}\phi= * F(A). \end{equation}
This is supplemented by asymptotic conditions at infinity. In the case when the structure group is $SU(2)$ one requires that $\vert\phi\vert$ tends to $1$ at infinity. Then the restriction of $\phi$ to a large sphere defines a map, up to homotopy,  from $S^{2}$ to $S^{2}$ and hence an integer degree $k$. 
For each $k\geq 1$ there is a moduli space $M_{k}$ of monopoles,  which has dimension $4k-1$. When $k=1$ the moduli space is just $\bR^{3}$ and the solution can be thought of as a particle centred at a point in $\bR^{3}$. One has a hyperk\"ahler picture, inherited from that in $\bR^{4}$, and the moduli spaces can be seen as hyperk\"ahler quotients. Again, this requires one to work with suitable framed moduli spaces $\tilde{M}_{k}$ which are circle bundles over the $M_{k}$, hence of dimension $4k$. There is also an analogue of the ADHM construction, due to Nahm, which gives another hyperkahler quotient description of the same moduli spaces.   

While the general set-up is very similar to the instanton case there are important differences. One is that the hyperk\"ahler metrics on the $M_{k}$ are complete, so the construction gives a new source of complete hyperk\"ahler manifolds.  This was developed in work of Atiyah and Hitchin \cite{kn:AH1}. The first interesting case is when $k=2$. The moduli space $\tilde{M}_{2}$ has dimension $8$ but another hyperk\"ahler quotient by the circle action gives a reduced space $\uM$ of dimension $4$ which can be identified, topologically, as the cotangent bundle of the real projective plane. The quotient construction involves fixing an origin $\bR^{3}$ and the solutions can be thought of, roughly speaking, as a pair of $1$-monopoles centred at antipodal points $\pm x\in \bR^{3}$ with a circle-valued  phase which determines how the $1$-monopoles are glued together. The line $\bR x$ gives a point in $\bR\bP^{2}$ and the phase and $\vert x\vert$ can be thought of as polar co-ordinates  on the cotangent space. This approximate description becomes more exact in the asymptotic regime where $\vert x\vert\rightarrow \infty$ and breaks down when $\vert x\vert\rightarrow 0$.

The group $SO(3)$ of rotations of $\bR^{3}$ acts by isometries on $\uM$. (It does not preserve the hyperk\"ahler structure $I,J,K$ but rotates it.)
The generic orbits have codimension $1$, so the metric can be regarded as a solution of a system of ODE's (the Einstein equations with this symmetry). Atiyah and Hitchin were able to solve this system in terms of elliptic integrals, so arriving at an explicit formula for the metric.  The Atiyah-Hitchin manifold $\uM$ has been  important in subsequent developments in 4-dimensional Riemannian geometry. The asymptotic behaviour is different from the ALE manifolds we discussed above: it is called asymptotically locally flat (ALF). Other, older,  examples  are the Taub-Nut space and families obtained by the Gibbons Hawking construction related to the $A_{k}$ series of ALE manifolds. The Atiyah-Hitchin manifold is the first in another $D_{k}$ series of ALF gravitational instantons which have been the scene for much subsequent activity, for example \cite{kn:Cherkis}.

In another direction, Atiyah and Hitchin used their description of the metric to investigate the \lq\lq dynamics'' of monopoles. To explain the idea, consider  a potential function $V$ on some Riemannian manifold $Q$ which has a minimum, say $V=V_{0}$, on a submanifold $L\subset Q$. We consider the motion of a particle $q(t)$ on $Q$ in the potential $V$, so the energy is
$$  V+ \frac{1}{2}\vert \dot{q}\vert^{2}. $$
If the energy is close to $V_{0}$ the particle is constrained to move in a small neighbourhood of $L$ and it can be shown that the motion is well-approximated by geodesic motion on $L$, for the induced Riemannian metric. In the case at hand one has an infinite-dimensional, field-theory, version of this set-up. The manifold $Q$ becomes the space of gauge equivalence classes of pairs $(A,\phi)$ on $\bR^{3}$ and the potential $V$ is
$$  V(A,\phi)= \int_{\bR^{3}} \vert F(A)\vert^{2} + \vert d_{A}\phi\vert^{2} $$
which is minimised on the monopole moduli space (similar to what we discussed for instantons). The equation of motion now becomes a hyperbolic equation on $\bR^{3,1}$. The principle above indicates that low energy solutions to this hyperbolic equation can be approximated by geodesics on the moduli space. Using their formulae for the metric and solving the geodesic equations numerically, Atiyah and Hitchin obtained a wealth of information about this dynamical question.

In other work, Atiyah initiated the study of monopoles on 3-dimensional hyperbolic space \cite{kn:MFA9}.
If we take a standard circle action on $S^{4}$, with fixed point set a $2$-sphere, the quotient  $(S^{4}\setminus S^{2})/S^{1}$ can naturally be viewed as hyperbolic 3-space $H_{3}$. More precisely,  $S^{4}\setminus S^{2}$ is conformally equivalent to $H_{3}\times S^{1}$. Just as translation-invariant instantons on $\bR^{4}$ correspond to monopoles on $\bR^{3}$ so is the same for $S^{1}$-invariant instantons and hyperbolic monopoles. The great advantage is that the analysis can be greatly simplified by considering instantons which extend smoothly over the 2-sphere. This is not the most general case: one has a real parameter given (for $SU(2)$ bundles) by the length of the Higgs field at infinity and solutions extend when this parameter is integral. For non-integral values one gets singular connections of the kind studied later by Kronheimer and Mrowka \cite{kn:KM}. Taking this parameter to infinity corresponds, under scaling,  to making the curvature of the hyperbolic space tend to $0$,  so the geometry converges to $\bR^{3}$ on bounded regions. Monopoles over both $\bR^{3}$ and $H_{3}$ have twistor descriptions. The relevant \lq\lq mini-twistor'' spaces, introduced by Hitchin \cite{kn:NJH1}, are the spaces of oriented geodesics and are complex surfaces. For $\bR^{3}$ this is the tangent bundle $TS^{2}$ of the 2-sphere and in the hyperbolic case it is the complement of the anti-diagonal $\overline{\Delta}$ in $S^{2}\times S^{2}$ (i.e. $\overline{\Delta}$ is the set of antipodal pairs). The tangent bundle $TS^{2}$ can be identified with the singular quadric cone in $\bC\bP^{3}$ minus the vertex while $S^{2}\times S^{2}$ is a smooth quadric,  so in 
the convergence of the hyperbolic to the flat case we encounter again the degeneration of smooth quadrics to the cone.

As a final example of the application of these hyperk\"ahler quotient ideas we mention Hitchin's work on Higgs bundles over Riemann surfaces. That is, pairs  $(E,\phi)$ where $E$ is a rank $r$ bundle and  $\phi\in H^{0}({\rm End } E\otimes K)$. The theory can also be seen as a dimension reduction of the instanton theory in four dimensions. The moduli space of Higgs bundles is hyperk\"ahler. It has another description as the moduli space of irreducible projective representations of the fundamental group in $GL(r,\bC)$, thus making a link with the old work of Weil \cite{kn:Weil}.

\subsubsection{Topology of instanton and monopole moduli spaces }

We move back a few years,  to a 1978 paper of Atiyah and Jones \cite{kn:AJ} and return to  the line of ideas in the paper of Atiyah and Bott. For a $G$-bundle $E$ over  $S^{4}$ we have the space of connections $\cA$ and gauge group $\cG$ and the Yang-Mils functional $\cL$ which takes its minimum on the instanton moduli space. The general question is whether the topology of $M$ can be related to that of $\cA/\cG$ and possible non-minimal critical points, as done by Atiyah and Bott in the 2-dimensional case. It will be convenient here to use the framed moduli spaces and hence restrict to the subgroup $\cG_{\infty}\subset \cG$ fixing the fibre over infinity. Then $\cG_{\infty}$ acts freely on $\cA$ and we do not have any trouble with reducible connections. In the 2-dimensional case we focused on rational cohomology but here the main interest is in more subtle, torsion,  phenomena. The equivalence class of the $G$-bundle $E$ is determined by an element of $\pi_{3}(G)$--the homotopy class of a transition function defined on an equatorial $3$-sphere. If, as we will assume, $G$ is a simple Lie group we have $\pi_{3}(G)\cong \bZ$ so the bundles are labelled by an integer $k$. 
The quotient $\cB=\cA/\cG_{\infty}$ is homotopy equivalent to the corresponding component  $\Omega^{3}_{k}(G)$ of the third loop space $\Omega^{3}G={\rm Maps}\ (S^{3},G)$,   where throughout we use based maps. This homotopy equivalence can be seen by considering preferred transition function defined by the  parallel transport of a connection along meridians.  Atiyah and Jones initiated the investigation of the maps on homotopy and homology induced by the inclusion
$\tM_{k}\rightarrow \Omega_{k}^{3}G$. (The group structure implies that all  the components $\Omega^{3}_{k}(G)$ have  the same homotopy type so we will sometimes blur the distinction between them.)

From now on we restrict to the group $G=SU(2)$, which is diffeomorphic to $S^{3}$.  The main result of Atiyah and Jones states that for each $q$, when $k$ is sufficiently large, the induced map $i_{*}:H_{q}(\tM_{k})\rightarrow H_{q}(\Omega^{3}_{k}(S^{3}))$ is a projection onto a direct summand. In particular, all the homology  of $\Omega_{k}^{3}(S^{3}$) is eventually seen in the instanton moduli spaces, as $k$ increases. Their proof used the family of  t'Hooft solutions to the instanton equations, mentioned in (2.1). These can be described in various ways. In terms of the ADHM data (30),(31) they arise when $r=2$ and the composite $P\circ Q$ is zero. So $\alpha_{1}, \alpha_{2}$ commute and for the t'Hooft solutions we take them to be diagonal 
$$\alpha_{1}={\rm Diag}(\lambda_{1}, \dots, \lambda_{k})\ \ \alpha_{2}={\rm Diag}(\mu_{1}, \dots, \mu_{k}). $$
The nondegeneracy conditions require that the $k$ points $(\lambda_{i}, \mu_{i})$ in $\bC^{2}$ are distinct and one finds that the remaining data is given by positive weights, or scales, associated to these $k$ points. For the topological discussion these scales can all be set to $1$ and we have a family of solutions parametrised by the configuration space $C_{k}(\bR^{4})$ of $k$ distinct points in $\bC^{2}=\bR^{4}$. Now Atiyah and Jones show that the map from $C_{k}(\bR^{4})$ to $\Omega^{3}_{k}(S^{3})$ given by this t'Hooft family is homotopy equivalent to a composite of two other maps. One is  the map
$$  \Omega^{4}(S^{4})\rightarrow \Omega^{4}(\bH\bP^{\infty})= \Omega^{4}BSU(2)=\Omega^{3}SU(2), $$
defined by the standard inclusion $S^{4}=\bH\bP^{1}\subset \bH\bP^{\infty}$. The other is the case $n=4$ of a map $C_{k}(\bR^{n})\rightarrow \Omega_{k}^{n}S^{n}$, well known in homotopy theory,  which Segal had shown  induces an isomorphism on $q$-dimensional homology  once $k>>q$. Combining this with other arguments from algebraic topology they deduced their result. 

The most influential part of the Atiyah-Jones paper was a conjecture (subsequently the \lq\lq Atiyah-Jones conjecture'') they made, to the effect that  the homology $H_{q}(\tM_{k})$ should be equal to $H_{k}(\Omega^{3}(S^{3}))$ once $k>>q$. Part of the motivation for this came from an important analogy between Yang-Mills theory in dimension $4$ and harmonic maps in dimensions $2$.  In this analogy the instantons correspond to holomorphic maps from a Riemann surface to a K\"ahler manifold. The simplest case is when both domain and target are the Riemann sphere, so we are considering rational maps which, after fixing base points, can be taken in the form
$$  f(z)= \sum_{i=1}^{k} \frac{a_{i}}{z-z_{i}}, $$ or limits of these when some of the points $z_{i}$ coincide. Segal showed  in \cite{kn:GBS}, around the same time as the Atiyah-Jones work, that the inclusion of the space $R_{k}$ of rational maps of degree $k$ into $\Omega^{2}_{k}S^{2}$ induces an isomorphism on $H_{q}$ once $k>>q$, which is the analogue of the  Atiyah-Jones conjecture.   

This conjecture of Atiyah and Jones was the scene for a lot of activity in the 1980's. It was eventually proved (for $SU(2)$ instantons) by Boyer, Hurtubise, Milgram and Mann in 1993 \cite{kn:BHMM}.  The activity followed two main lines: variational methods, chiefly by Taubes, and more explicit geometrical descriptions of the moduli spaces.  We recall some background for the variational problem. The Yang-Mills functional is defined in all dimensions $n$. Like many other variational problems, there is a a critical dimension $n=4$ where it is scale invariant. For harmonic maps the corresponding dimension is $2$. Below this critical dimension the analysis is relatively straightforward. So, for example, Rade established long-time existence and convergence of the Yang-Mills flow (24) in dimensions 2 and 3.  Roughly speaking, in dimensions 2 or 3 the topology of the space of connections modulo gauge must match up with the critical points in the usual Morse theory fashion. 

Going to dimension 4, with an eye to the Atiyah-Jones conjecture, there is a network of difficulties. Most important, at the time Atiyah and Jones were writing there was no theory which could implement variational arguments at the critical dimension. The basic problem is the lack of compactness: as sequence of connections with bounded Yang-Mills functional may fail to converge due to \lq\lq bubbling'' at points. In the following years such results were developed,  using deeper analysis of  this bubbling. In the case of harmonic maps of surfaces this development was initiated by Sacks and Uhlenbeck and in the Yang-Mills case by Taubes and Uhlenbeck  The second difficulty, at the time Atiyah and Jones were writing, was that little was known about the higher critical points of the Yang-Mills functional. The analogy with harmonic maps gave some grounds for thinking that there might be none. Subsequently these were proved to exist, first by Uhlenbeck, Sibner and Sibner \cite{kn:USS} using variational arguments and later by relatively  explicit constructions \cite{kn:Sadun}. But an important result of Taubes \cite{kn:CHT1} showed that the index of these higher  critical points must increase with $k$: the index cannot be less than roughly one quarter the dimension of the moduli space of instantons. Taubes' argument used the action of quaternions in a way remininscent of the Morse theory proof of the Lefschetz hyperplane theorem. This index bound is consistent with the Atiyah-Jones conjecture and Taubes proved a number of results in that direction \cite{kn:CHT2}, but not  the full conjecture. More precisely, Taubes introduced maps $T_{k}:\tM_{k}\rightarrow \tM_{k+1}$, unique  up to homotopy, and showed that resulting direct limit of homology groups is $H_{*}(\Omega_{0}^{3}S^{3})$. In particular, {\it if} the homology groups $H_{q}(\tM_{k})$ stabilise as $k\rightarrow \infty$ then the limit must be $H_{q}(\Omega_{0}^{3}S^{3})$. But this leaves open the possibility that there are additional classes in $H_{q}(\tM_{k})$, for arbitrarily large $k$, which map to zero in $H_{q}(\tM_{k'})$ for some $k'>k$.  

We turn now to the more \lq\lq geometrical'' approaches. Gravesen \cite{kn:Gravesen} obtained a result similar to that discussed above of Taubes by using a description due to Atiyah \cite{kn:MFA8} of the $G$-instanton moduli spaces as spaces of holomorphic maps from $S^{2}$ to the loop space $\Omega G$. This is an infinite dimensional complex manifold and Gravesen's work extends  Segal's theory for rational maps to the infinite dimensional situation. The technique of Boyer{\it et al} in \cite{kn:BHMM} exploited the identification of instantons with holomorphic bundles  on the plane mentioned in 2.4.3 above. They defined a stratification of the moduli space based on the configuration of jumping lines and, by a detailed study of this, were able to show that the homology groups stabilise. Combined with the results of Taubes or Gravesen, this gave  a proof of the Atiyah-Jones conjecture. Another proof, of a slightly different result, was found by Kirwan using the ADHM description and her general theory of quotient spaces \cite{kn:FCK3}.

 The topology of instanton moduli spaces and the variational theory of the Yang-Mills functional seems    an interesting area for future study. Do  analogues of  the Atiyah-Jones conjecture hold true for other 4-manifolds? One recent relevant  development is due to Waldron, who showed that the Yang-Mills flow in dimension $4$ has long-time smooth solutions \cite{kn:AW}.

While variational methods  have, so far, been only partially successful for Yang-Mills theory on 4-manifolds, there is a very satisfactory  picture for analogous questions involving monopoles.  Taking the gauge group $SU(2)$ the moduli space of charge $k$ monopoles can be identified with the rational map spaces $R_{k}$. Taubes' index bounds apply also to monopoles and he developed a full variational theory in this case \cite{kn:CHT2}. Putting this together gives an alternative proof of Segal's result on the homology approximation of $\Omega^{2}_{k}S^{2}$ by $R_{k}$ for large $k$. Similar results hold for other structure groups where the  moduli spaces that appear are spaces of holomorphic maps from $S^{2}$ to coadjoint orbits. Monopoles also tie up with the instanton discussion. The proof, mentioned above,  by Uhlenbeck, Sibner and Sibner of the existence of higher critical points can be seen as a variant, for hyperbolic monopoles, of Taubes' theory for Euclidean monopoles.

\subsubsection{Localisation}

Atiyah was very interested in an integration formula found by Duistermaat and Heckman in 1982 \cite{kn:DH}. Consider a compact sympletic $2m$-manifold $(V,\omega)$ with a circle action generated by a Hamiltonian function $H$. Suppose for simplicity that the  action has a finite set of fixed points, or in other words $H$ has non-degenerate critical points. The Duistermaat-Heckman formula is
\begin{equation}  \int_{V} e^{-itH} \frac{\omega^{2m}}{m!} = \sum_{p} \frac{e^{-itH(p)}}{e(p) t^{m}}, \end{equation}
where $p$ runs over the fixed points and $e(p)$ is the product of the weights of the action on $TV_{p}$. 

If we take any function $H$ on $V$, standard  theory gives an {\it asymptotic} formula as $t\rightarrow \infty$ for the integral on the left hand side of (33) with a sum of contributions from the critical points which  reproduces the expression on the right hand side in this case. This is the stationary phase approximation.  The point of the Duistermaat-Heckman result is that the fact that $H$ generates a circle action implies that the stationary phase approximation is exact.

Considering $H$ as a map $H:V\rightarrow \bR$, let $\mu$ be the measure on $\bR$ given by the push forward of the volume form $\omega^{m}/m!$ on $V$. The integral on the left hand side of (33) now appears as the Fourier transform of $\mu$ and simple arguments from analysis, together with the asymptotic statement, show that the equation (33) is equivalent to the fact that $\mu$ is piecewise polynomial. In the simplest case, when $V$ is the $2$-sphere the measure $\mu$ is supported in an interval on which it is a constant multiple of Lebesgue measure. This is a restatement of Archimedes' theorem on the area of zones in the sphere. The same holds for the action of 
of a torus $T^{k}$ on a symplectic manifold $(V,\omega)$ with a  moment map $m:V\rightarrow \bR^{k}$, where $\bR^{k}$ is viewed as the dual of the Lie algebra. The push-forward of the volume form is piecewise polynomial. Using ideas from Morse Theory, Atiyah obtained a result in similar vein: the support of this measure---that is,  the image of $m$---is a convex polyhedron \cite{kn:MFA6}. (Similar results were obtained about the same time by Guillemin and Sternberg.) In the case when $V$ is a co-adjoint orbit of a compact Lie group and $T^{k}$ is the maximal torus this recovers  classical matrix inequalities and has an intimate connection with representation theory.

In another direction, Atiyah and Bott explained how the formula (33) can be regarded as a localisation formula in equivariant cohomology \cite{kn:AB2}. Recall that for any  group $G$ acting on $V$ the equivariant cohomology $H^{*}_{G}(V)$ is the cohomology of the space $V_{G}=EG\times_{G}V$ which fibres over $BG$ with fibre $V$. For a compact oriented $2m$ manifold $V$ we have then an integration-over-the-fibre map
$$  I: H^{*}(V_{G})\rightarrow H^{*-2m}(BG). $$
In the case at hand, with $G=S^{1}$ the cohomology $H^{*}(BG)$ is the polynomial ring in a generator $\tau$ of degree $2$ so, taking complex co-efficients,  we have a map
$$  H^{*}_{S^{1}}(V)\rightarrow \bC[\tau]. $$
Atiyah and Bott showed that, in the abstract, this map can be computed from data at the fixed points of the action. The simplest case is when the action is free. Then the equivariant cohomology is the same as the cohomology of the quotient space $V/S^{1}$ and it is clear that the integration over the fibre vanishes.  The general Atiyah-Bott formula, for isolated fixed points, is 
\begin{equation}  I(\alpha)= \sum_{p}  \frac{i_{p}^{*}(\alpha)}{E(TV_{p})} \end{equation} where $i_{p}$ is the inclusion of $p$ in $V$ and $E(TV_{p})$ is the equivariant Euler class of the tangent bundle. Each term in the sum on the right hand side is computed as an element of the field of fractions $\bC(\tau)$ but the sum lands in the polynomials $\bC[\tau]$.

To connect with the Duistermaat-Heckman formula we recall the Cartan model for equivariant cohomology via differential forms. In the case of a circle action we take $\Omega^{*}_{V}\otimes \bC[\tau]$ with a differential
$$ d_{X}=  d+ \tau i_{X}, $$
where $i_{X}$ is contraction with the vector field $X$ generating the action. Acting on the $S^{1}$-invariants in  $\Omega^{*}_{V}\otimes \bC[\tau]$, this has square zero and the resulting cohomology groups are isomorphic to $H^{*}_{S^{1}}(V)$. Now the hypothesis that $H$ is a Hamiltonian for the circle action is equivalent to saying that $\omega+ \tau H$ is a $d_{X}$-closed form. In other words $H$ gives an extension of $\omega$, which is an ordinary closed form, to an equivariantly closed form. Now the integral on the left hand side of (33) is just 
$$  I(\exp(\omega+ \tau H))\in \bC[\tau], $$
except that we replace the formal variable $\tau$ by the complex variable $it$ in interpreting the formula. 

\subsubsection{The Seiberg-Witten curve and Nekrasov's instanton counting} 

We finish this article with a brief discussion of a development which brings together many of the themes we have discussed in Atiyah's work.

As we mentioned above, the individual terms in the Atiyah-Bott fixed point formula are defined as rational functions of $\tau$. The same applies for an action of a torus $T=T^{\mu}$ with isolated fixed points, where each fixed point contributes  a rational function of variables $(\tau_{1}, \dots, \tau_{\mu})$. The fixed point contribution
$$  \frac{i_{p}^{*}(\alpha)}{E(TV_{p})}$$
makes sense for classes $\alpha\in H^{*}_{T}(V)$ which are of the wrong dimension to integrate over the manifold. In particular we can consider the class $1$ in $H^{0}$. For a compact manifold $V$ the sum of these contributions must vanish but we will now consider the case of an action on a non-compact manifold with a finite number of fixed points. Then we can make the formal definition
   \begin{equation}  \int_{V} 1 = \sum_{p} \frac{1}{E(TV_{p})} \end{equation}
which is a rational function of $(\tau_{1},\dots,\tau_{\mu})$. 

To understand this more, suppose that $V$ is a manifold with boundary $\partial V=W$ on which the torus $T$ acts and there are no fixed points on the boundary. Then the quantity in (35) is a topological invariant of the action of $T$ on  the boundary $W$. For the simplest case, suppose that $T=S^{1}$ and that the action on $W$ is free. Then $W/S^{1}=Q$ is a $(2m-2)$-manifold and $W\rightarrow Q$ is an $S^{1}$-bundle with a Chern class $c_{1}\in H^{2}(Q)$. Then it is not hard to see that

\begin{equation}  \int_{V} 1= \tau^{-m} \langle c_{1}^{m-1}, [Q]\rangle \end{equation}

In the general case of a $T$-action with no fixed points on $W$, for almost all elements $\xi$ of the Lie algebra of $T$ the corresponding vector field $X_{\xi}$ on $W$ has no zeros. Let $A$ be a $1$-form on $W$ such that the Lie derivative of $A$ under $X_{\xi}$ vanishes and $A(X_{\xi})=1$. Such $1$-forms exist: for example we can take a $T$-invariant Riemannian metric on $W$ such that $X_{\xi}$ has length $1$ and the $1$-form dual to $X_{\xi}$. We then obtain a quantity
$$  J= \int_{W} A\wedge dA^{m-1}, $$
and it is straightforward to show that this is independent of the choice of $A$ and is a rational function of $\xi$. Then one has a formula
\begin{equation}  \int_{V} 1= J \end{equation}
where the right hand side is a rational function defined  on the Lie algebra of the torus. (In the case when $W$ is 3-dimensional, this invariant is related to the asymptotic Hopf invariant introduced by Arnold \cite{kn:Arn}.)

Nekrasov applied these ideas in \cite{kn:NN} to 
framed moduli spaces  $\tilde{M}_{k,r}$ of instantons on $\bR^{4}$ with rank $r$ and Chern class $k$. There are commuting actions of $SO(4)$, acting by rotations of $\bR^{4}$,  and of $U(r)$, acting on the framing. We take  maximal tori $T^{2}$ in $SO(4)$ and  $T^{r-1}\in SU(r)\subset U(r)$ (the centre of $U(r)$ acts trivially). More precisely, Nekrasov considered a \lq\lq completion'' $\oM_{k,r}$ of the moduli space which depends on a choice of complex structure $\bR^{4}=\bC^{2}$. As we mentioned in 2.4.3 the framed instanton moduli spaces can be identified with moduli spaces of holomorphic bundles on the complex plane  trivialised at infinity. Then the completion  $\oM_{k,r}$ can be defined as a set of Gieseker stable torsion-free sheaves on the plane, framed at infinity. From another point of view, we modify the ADHM equations (30,(31) to
$$   [\alpha_{1}, \alpha_{2}]+PQ=0\ \ \ , \ \ \ [\alpha^{*}_{1},\alpha_{1}]+ [\alpha^{*}_{2},\alpha_{2}]+ P^{*}P-QQ^{*}=\epsilon, $$
for some nonzero $\epsilon$. In other words we are changing the level set in the definition of the hyperk\"ahler quotient, just as we have discussed before. In the simplest case of $k=1, r=2$ the moduli space is 
$$  M_{2,1}=\bC^{2}\times (\bC^{2}\setminus \{0\})/\pm 1) . $$ 
The first factor gives the \lq\lq centre'' of the instanton and the norm of the second factor gives the \lq\lq scale''. The \lq\lq Uhlenbeck completion'', which does not depend on a choice of complex structure, adds ideal instantons with scale zero to give
$$  \bC^{2}\times (\bC^{2}/\pm 1), $$
whereas the Gieseker completion is obtained by blowing up the singular set to add an exceptional divisor $\bC^{2} \times \bC\bP^{1}$.

The advantage of the Gieseker completion is that the spaces $\oM_{r,k}$ are smooth manifolds. The action of $T^{2}\times T^{r-1}$ on the moduli space extends to $\oM_{r,k}$ and Nekrasov showed that there are a finite number of fixed points which can be described explicitly in terms of Young diagrams. So the preceding discussion applies. For simplicity we now take $r=2$ so the torus acting is $T^{2}\times S^{1}$ and we take co-ordinates $(\epsilon_{1}, \epsilon_{2}, a)$ on the Lie algebra.  Nekrasov defines
\begin{equation}   Z(\epsilon_{1}, \epsilon_{2}, a, \Lambda)= \sum_{k=0}^{\infty} \Lambda^{k} \int_{\oM_{2,k}} 1. \end{equation}
 These definitions were motivated by supersymmetric Yang-Mills theory and in particular the work of Seiberg and Witten.  Nekrasov made a series of  conjectures about this function $Z$. One was that  $$
F(\epsilon_{1},\epsilon_{2},a,\Lambda)= \epsilon_{1}\epsilon_{2}\log Z(\epsilon_{1},\epsilon_{2},a)$$
is  holomorphic across $\epsilon_{1}, \epsilon_{2}=0$.  Another was that $F(0,0,a,\Lambda)$ is given by a certain formula involving the periods of the meromorphic 1-form $ zdw/w$ on the \lq\lq Seiberg-Witten curve'': the elliptic curve $C_{u,\Lambda}$ defined by the solutions $(z,w)$ of the equation
\begin{equation} \Lambda (w+w^{-1}) = z^{2} + u . \end{equation}

These conjectures were proved subsequently by Nekrasov and Okounov \cite{kn:ON} and by Nakajima and Yoshioka  \cite{kn:NY} and we refer to those papers for further details. For the present article, the point is that this work---combining Atiyah-Bott localisation, instanton moduli spaces, the ADHM constuction and hyperk\"ahler quotients---produced fundamental results in quantum Yang-Mills theory.




\begin{thebibliography}{99}
\bibitem{kn:Arn}  V. I. Arnold {\em The asymptotic Hopf invariant and its applications}\  Proc. Summer School in
 Diff. Equations at Dilizhan, 1973; English translation  Sel. Math.
 Sov.  327-345 (1986)
 \bibitem{kn:ADK} { A. Asok, B. Doran and F. Kirwan} \ {\em Yang-Mills theory and Tamagawa numbers}\  Bull. Lond. Math. Soc. 40 533-67 (2008)
\bibitem{kn:MFA1}{ M.F Atiyah}\ {\em Complex fibre bundles and ruled surfaces}\  Proc. Lon. Math. Soc. I 407-34 (1955)
\bibitem{kn:MFA2}{ M.F Atiyah}\  {\em Complex analytic connections in fibre bundles}\  Trans. Amer. Math. Soc. 85 181-207 (1957)
\bibitem{kn:MFA3}{ M.F Atiyah}\ {\em  Vector bundles over an elliptic curve}\  Proc. Lond. Math. Soc. VII 414-52 (1957) 
\bibitem{kn:MFA4}{ M.F Atiyah} \ {\em On analytic surfaces with double points} Proc. Roy. Soc. A 247 237-44 (1958)
\bibitem{kn:MFA5}{ M.F Atiyah} \ {\em Geometry of Yang-Mills fields}\  Lezioni Fermiane Acad. Nazionale dei Lincei \& Scuola Normale Sup. Pisa (1979)
\bibitem{kn:MFA7}{ M.F Atiyah}\ {\em  Green's functions for self-dual four-manifolds}
Advances in Mathematics Supplementary studies 6A 129-58 (1981)
\bibitem{kn:MFA6}{ M.F Atiyah} \ {\em Convexity and commuting Hamiltonians}\  Bull.
Lon. Math. Soc. 14 1-15 (1982)
\bibitem{kn:MFA8}{ M.F. Atiyah} \ {\em  Instantons in two and four dimensions}\  Commun.Math.Phys. 93 437-51 (1984)
\bibitem{kn:MFA9}{ M.F Atiyah}\ {\em  Magnetic monopoles in hyperbolic space}\  in Proc. Bombay Colloquium on vector bundles over algebraic varieties  Oxford UP 1-34 (1987)
\bibitem{kn:AB1}{ M.F Atiyah and R. Bott} \ {\em The Yang-Mills equations over Riemann surfaces}\  Phil. Trans. R. Soc. Lon. A 308 523-615 (1982)
\bibitem{kn:AB2}{ M.F Atiyah and R. Bott}\ {\em  The moment map and equivariant cohomology}\  Topology 23 1-28 (1984) 
\bibitem{kn:ADHM}{ M.F Atiyah, V. G. Drinfeld, N. J. Hitchin and Yu. I Manin }\ {\em  Construction of instantons} \  Phys. Lett. 65A 185-7 (1978)
\bibitem{kn:AH1}{ M.F Atiyah and N. J. Hitchin}\ {\em  Geometry and dynamics of magnetic monopoles}\  Princeton U. P. (1988) 
\bibitem{kn:AHS1}{ M.F Atiyah, N. J. Hitchin and I. M. Singer}\ {\em  Deformations of instantons}\  Proc. Nat. Acad. Sci USA 74 2662-3 (1977)
\bibitem{kn:AHS2}{ M.F Atiyah, N. J. Hitchin and I. M. Singer}\ {\em  Self-duality in four-dimensional Riemannian geometry}\  Proc. Roy. Soc. Lond. A 362 425-61 (1978)
\bibitem{kn:AJ} { M. F. Atiyah and J. D. S. Jones} {\em Topological aspects of Yang-Mills theory}\  Commun. Math. Phys. 61 97-118 (1978)
\bibitem{kn:MFAW} { M. F. Atiyah and R. S. Ward}\ {\em  Instantons and algebraic geometry} Commun. Math. Phys. 55 117-24 (1977)
\bibitem{kn:Bateman} { H. Bateman} \ {\em The solution of partial differential equations by means of definite integrals}\  Proc. Lond. Math. Soc. 1 451-458 (1904)
\bibitem{kn:BH} { W. Barth and K. Hulek}\ {\em   Monads and moduli of vector bundles}\  Manuscripta Math. 25  323-347 (1978)
\bibitem{kn:BHMM} { C .P. Boyer, J.C. Hurtubise, B.M. Mann and R.J. Milgram} {\em The topology
of instanton moduli spaces. I. The Atiyah-Jones conjecture.}\  Ann. of Math.
 137  561-609 (1993)
 \bibitem{kn:Cherkis}  { S. Cherkis and A. Kapustin} {\em $D_k$ gravitational instantons and Nahm equations}\  Adv. Theor. Math. Phys. 2  1287-1306 (1999).
\bibitem{kn:SKD1}{ S.K.Donaldson} {\em  Instantons and geometric invariant theory}\  Comm. Math. Phys. 93  453-460 (1984)
\bibitem{kn:SKD2}{S. K. Donaldson} {\em Nahm's equations and the classification of monopoles}\  Commun. Math. Phys. 96 387-407. (1984)
\bibitem{kn:SKD3} { S.K. Donaldson} {\em Topological field theories and formulae of Casson and Meng-Taubes}\   Geom. Topol. Monogr. 2 87-102 (1999)
\bibitem{kn:SKD4}{ S.K. Donaldson}\ {\em  Moment maps and diffeomorphisms}
Sir Michael Atiyah: a great mathematician of the twentieth century.
Asian J. Math. 3 1-15 (1999)
\bibitem{kn:SKD5} {\em S. K. Donaldson}\  {\em Moment maps in differential geometry}\ 
Surveys in differential geometry, Vol. 8 171-189 (2002) 
\bibitem{kn:DH} { J. J. Duistermaat and G.J. Heckman}{\em  On the variation in the cohomology of the symplectic form of the reduced phase space}\  Invent. Math. 69 (1982)
\bibitem{kn:EK} { R. Earl and F.C.Kirwan} \ {\em  Complete
sets of relations in the cohomology rings of moduli spaces of holomorphic
bundles and parabolic bundles over a Riemann surface}\  Proc. London Math.
Soc. 89  570-62  (2004)
\bibitem{kn:GL} { D. Gaitsgory and J. Lurie} {\em
Weil's conjecture for function fields. Vol. 1.}
Annals of Mathematics Studies, 199. Princeton University Press (2019)
\bibitem{kn:WMG} { W. M. Goldman} {\em  The symplectic nature of fundamental groups of surfaces}\  Adv. in Math. 54  200-225 (1984)
\bibitem{kn:Gravesen} { J. Gravesen}
{\em On the topology of spaces of holomorphic maps}\ 
Acta Math. 162  247-286 (1989)
\bibitem{kn:GH} { P. Griffiths and J. Harris} \ {\em Principles of Algebraic Geometry}\   Wiley 1978
\bibitem{kn:HN} { G. Harder and M.S. Narasimhan}
\ {\em On the cohomology groups of moduli spaces of vector bundles on curves}
Math. Ann. 212  215-248 (1974/75)
\bibitem{kn:NJH1} { N. J. Hitchin} {\em Monopoles and geodesics}\  Comm. Math. Phys. 83  579-602 (1982)
\bibitem{kn:NJH2} { N. J. Hitchin} {\em  The self-duality equations on a Riemann surface}\  Proc. London Math. Soc.  55  59-126  (1987)
\bibitem{kn:HKLR} { N. J. Hitchin, A. Karlhede, U. Lindstr\"om, M. Rocek}
\ {\em Hyper-K\"ahler metrics and supersymmetry}\  Comm. Math. Phys. 108  535-589
(1987)
\bibitem{kn:Horrocks}  { G. Horrocks}\ {\em  Vector bundles on the punctured spectrum of a local ring} Proc. London Math. Soc.  14  689-713 (1964)
\bibitem{kn:JK}{L. Jeffrey and F. C. Kirwan}\ {\em  Intersection theory on moduli
spaces of holomorphic bundles of arbitrary rank on a Riemann surface}\  Annals
of Math. 148 109-96 (1998)
\bibitem{kn:KN} { A. D. King  and P. E. Newstead} \ {\em On the cohomology of the moduli space of rank 2 vector bundles on a curve}\   Topology 37 (1998) 407-418 (1998)
\bibitem{kn:FCK1}{F. C. Kirwan} \ {\em Cohomology of quotients in symplectic and algebraic geometry}\  Math. Notes 31 Princeton U.P. (1984)
\bibitem{kn:FCK2}{ F. C. Kirwan}\ {\em
The cohomology rings of moduli spaces of bundles over Riemann surfaces}\
J. Amer. Math. Soc. 5  853-906(1992)
\bibitem{kn:FCK3} {F.C.  Kirwan} \ {\em Geometric invariant theory and the Atiyah-Jones
conjecture}\  The Sophus Lie Memorial Conference (Oslo, 1992), 161-186, Scand.
Univ. Press, Oslo, 1994. 570-622.
\bibitem{kn:Kronheimer} P. B. Kronheimer \ {\em  The construction of ALE spaces as hyper-k\"ahler quotients}\  Jour. Diff. Geom. 29 665-683 (1989)
\bibitem{kn:KM} { P. B. Kronheimer and T. S. Mrowka}
{\em Gauge theory for embedded surfaces. I}\ 
Topology 32  773-826 (1993)
\bibitem{kn:Mumford} { D. B. Mumford}\ {\em  Projective invariants of projective structures and applications}\  1963 Proc. Internat. Congr. Mathematicians (Stockholm, 1962) pp. 526-530 Inst. Mittag-Leffler, Djursholm
\bibitem{kn:NY} {H. Nakajima and K. Yoshioka} \ {\em Instanton counting on blow-up
I: 4-dimensional pure gauge theory}\   Inventiones Math. 162 313-355 (2005)
\bibitem{kn:MSN} { M.S. Narasimhan}\ {\em   Geometry of moduli spaces of vector
bundles}\  Actes du Congr\`es International des Math\'ematiciens (Nice, 1970),
Tome 2,  Gauthier-Villars, Paris, 1971 (199-201)
\bibitem{kn:MSNR} { M.S. Narasimhan and S. Ramanan}\ {\em  Moduli of vector bundles
on a compact Riemann surface}\  Ann. of Math.  89  14-51  (1969)
\bibitem{kn:MSNCS} { M. S. Narasimhan and C. S. Seshadri}\ {\em  Stable and unitary vector bundles over an algebraic curve}\  Ann. of Math. 82 (1965) 540-567.
\bibitem{kn:NN}{ N. Nekrasov}\ {\em  Seiberg-Witten prepotential from instanton counting}\  Advances Theor. Math. Physics 7 831-64 (2003)
\bibitem{kn:ON}{ A. Okounkov and N. Nekrasov}\ {\em  Seiberg-Witten prepotential
and random partitions} In:
The unity of mathematics, 525-596, Progr. Math., 244, Birkh\"auser  (2006)
\bibitem{kn:Newstead1} { P. E. Newstead}\ {\em  Topological properties of some spaces of stable bundles} Topology 6 241-62 (1967)
\bibitem{kn:Newstead2} { P. E. Newstead}\ {\em  Stable bundles of rank two and odd degree over a curve of genus two}\  Topology 7 205-215 (1968)
\bibitem{kn:Rade} { J. Rade} \ {\em  On the Yang-Mills heat equation in two and three dimensions}\  J. Reine Angew. Math. 431 123-63 (1992)
\bibitem{kn:Sadun} { L. Sadun and J. Segert}
{\em Non-self-dual Yang-Mills connections with nonzero Chern number}\
Bull. Amer. Math. Soc.  24  163-70 (1991)
\bibitem{kn:SMS1} { S. M. Salamon} 
{\em Quaternionic K\"ahler manifolds}\ 
Invent. Math. 67 143-171 (1982)
\bibitem{kn:SMS2} S. M. Salamon \ {\em Differential geometry of quaternionic manifolds}\ Ann. Sci. \'Ecole Norm. Sup.  19  31-55 (1986) 
\bibitem{kn:GBS}  { G. B. Segal}\ {\em  The topology of spaces of rational functions}\  Acta Math. 143  39-72 (1979)
\bibitem{kn:RLES} { R. L. E. Schwarzenberger}\ {\em  Vector bundles on the projective plane}\  Proc. Lond. Math. Soc. 311 623-40 (1961)
\bibitem{kn:CHT1} { C.H. Taubes}\ {\em  Stability in Yang-Mills theories}
Comm. Math. Phys. 91 (1983) 235-263  (1983)
\bibitem{kn:CHT2} { C. H. Taubes} \ {\em Min-max theory for the Yang-Mills-Higgs
equations}\  Comm. Math. Phys. 97 473-540(1985)
\bibitem{kn:CHT3} { C. H. Taubes}\ {\em  The stable topology of self-dual moduli spaces}\  J. Differential Geom. 29  163-230  (1989) 
\bibitem{kn:MT} { M. Thaddeus} \ {\em Conformal field theory and the cohomology of the moduli space of stable bundles}\  J. Differential Geom. 35 (1992) 131-149(1992)
\bibitem{kn:MT2} { M. Thaddeus} {\em 
 Stable pairs, linear systems and the Verlinde formula}\  Invent. Math. 117 317-353 (1994)
\bibitem{kn:USS} { K. K. Uhlenbeck, L. Sibner and R. Sibner} \ {\em  Solutions to Yang-Mills equations that are not self-dual}\  Proc. Nat. Acad. Sci. U.S.A. 86 8610-8613 (1989)
\bibitem{kn:AW} { A. Waldron} \ {\em Long-time existence for Yang-Mills flow}
Invent. Math. 217  1069-1147 (2019)
\bibitem{kn:Ward} { R. S. Ward} \ {\em   On self-dual gauge fields}\  Phys. Lett. A 61 81-82 (1977)
\bibitem{kn:Weil} { A. Weil}\ {\em  Generalisation des fonctions abeliennes}\  J. Math. Pures Appl. 17 47-87 (1938)
\bibitem{kn:Witten} { E. Witten} {\em  On quantum gauge theories in two dimensions}
Comm. Math. Phys. 141 153-209  (1991)
\bibitem{kn:SW}{ S. Wolpert} {\em On the symplectic geometry of deformations of a hyperbolic surface}\ 
Ann. of Math.  117  207-234(1983)
\end{thebibliography}
\end{document}